\documentclass[12pt,reqno]{amsart}
\usepackage{amsmath,amssymb,amsthm}
\usepackage{amsfonts,amsxtra}
\usepackage{psfrag,graphicx}
\usepackage{hyperref}

\DeclareFontFamily{OML}{rsfs}{\skewchar\font'177}
\DeclareFontShape{OML}{rsfs}{m}{n}{ <5> <6> rsfs5 <7> <8> <9> rsfs7
  <10> <10.95> <12> <14.4> <17.28> <20.74> <24.88> rsfs10 }{}
\DeclareMathAlphabet{\mathfs}{OML}{rsfs}{m}{n}

\newcommand{\Z}{\ensuremath{\mathbb{Z}}}

\newcommand{\Q}{\ensuremath{\mathbb{Q}}}
\newcommand{\R}{\ensuremath{\mathbb{R}}}

\newcommand{\uDS}{\ensuremath{\mathcal{X}}}
\newcommand{\DS}{\ensuremath{\mathcal{D}}}

\newcommand{\iv}{\ensuremath{o}}
\newcommand{\tv}{\ensuremath{t}}

\newcommand{\As}{{\mathfs A}}
\newcommand{\co}{\colon\thinspace}

\newcommand{\la}{\langle}
\newcommand{\ra}{\rangle}

\theoremstyle{plain}
\newtheorem{theorem}{Theorem}
\newtheorem{lemma}[theorem]{Lemma}
\newtheorem{proposition}[theorem]{Proposition}
\newtheorem{corollary}[theorem]{Corollary}

\theoremstyle{definition}
\newtheorem{definition}[theorem]{Definition}

\newtheorem{example}[theorem]{Example}

\newtheorem{remark}[theorem]{Remark}
\newtheorem{question}[theorem]{Question}
\newtheorem*{ack}{Acknowledgements}

\numberwithin{theorem}{section}

\DeclareMathOperator{\Aut}{Aut}
\DeclareMathOperator{\Inn}{Inn}
\DeclareMathOperator{\Out}{Out}
\DeclareMathOperator{\SOut}{\Sigma Out}

\DeclareMathOperator{\red}{red}
\DeclareMathOperator{\col}{col}

\begin{document}

%%%%%%%%%%%%%%%%%%%%%%%%%%%%%%%%%%%%%%%%%%%%%%%%%%%%%%%%%%%%%%%%%%%%%%%%%%%%% 

\title[Automorphisms of Baumslag--Solitar groups]{Deformation spaces
  of $G$--trees and automorphisms of Baumslag--Solitar groups}
\author{Matt Clay} \date{\today}

\address{Department of Mathematics\\ University of Oklahoma\\ Norman,
OK\\ 73019-0315\\ USA}

\email{mclay@math.ou.edu}

\begin{abstract}
  We construct an invariant deformation retract of a deformation space
  of $G$--trees.  We show that this complex is finite dimensional in
  certain cases and provide an example that is not finite dimensional.
  Using this complex we compute the automorphism group of the
  classical non-solvable Baumslag--Solitar groups $BS(p,q)$.  The most
  interesting case is when $p$ properly divides $q$.  Collins and
  Levin computed a presentation for $\Aut(BS(p,q))$ in this case using
  algebraic methods.  Our computation uses Bass--Serre theory to
  derive these presentations.  Additionally, we provide a geometric
  argument showing $\Out(BS(p,q))$ is not finitely generated when $p$
  properly divides $q$.
\end{abstract}

\maketitle

%%%%%%%%%%%%%%%%%%%%%%%%%%%%%%%%%%%%%%%%%%%%%%%%%%%%%%%%%%%%%%%%%%%%%%%%%%%%% 

\section*{Introduction}

Baumslag--Solitar groups have the following standard presentations:
\begin{equation}\label{eq:bs}
BS(p,q) = \la x,t \ | \ tx^pt^{-1} = x^q \ra.
\end{equation}
When $p$ properly divides $q$ there are infinitely many similar
presentations for $BS(p,q)$ which highlights additional symmetries.
These groups were first studied by Baumslag and Solitar as some
examples of non-Hopfian groups \cite{ar:BS62}.  Our interest is in the
automorphism and outer automorphism groups, $\Aut(BS(p,q))$ and
$\Out(BS(p,q))$ respectively, of the non-solvable Baumslag--Solitar
groups.  For non-solvable Baumslag--Solitar groups neither $|p|$ nor
$|q|$ equals 1.  By interchanging $t \leftrightarrow t^{-1}$, we can
always assume that $|q| \geq p > 1$.

Presentations for these automorphism groups are known.  The first
result was by Collins, who gave a finite presentation for
$\Aut(BS(p,q))$ when $p$ and $q$ are relatively prime \cite{ar:C78}.
This result was extended by Gilbert, Howie, Metaftsis and Raptis to
the cases when $p$ does not divide $q$ or when $p = |q|$
\cite{ar:GHMR00}.  Collins and Levin had earlier studied the most
interesting cases, which is when when $p$ properly divides $q$
\cite{ar:CL83}.  In this case $\Aut(BS(p,q))$ is not finitely
generated.  A summary of these results appears in Section
\ref{sc:dbspq}.  Although we do not consider the solvable case, we
note that Collins found a finite presentation for the automorphism
group $\Aut(BS(1,q))$ that depends on the prime factorization of $q$
\cite{ar:C78}.

One of the purposes of this paper is to a give unified approach to the
computation of these automorphism groups.  To this end, we construct a
tree $X_{p,q}$ on which $\Out(BS(p,q))$ acts.  In the cases Gilbert et
al.~considered, the tree $X_{p,q}$ is a single point.  In the more
interesting cases, when $p$ properly divides $q$, this tree is
nontrivial.  The arguments from Gilbert et al.~are used to compute the
vertex stabilizers of the tree $X_{p,q}$.  Using Bass--Serre theory in
the case where $p$ properly divides $q$ we recover the presentations
of $\Aut(BS(p,q))$ and $\Out(BS(p,q))$ (Theorem
\ref{th:presentations}) originally found by Collins and Levin
\cite{ar:CL83}.  Prior to finding these presentations, we present a
simple geometric argument showing that $\Out(BS(p,q))$ (and hence
$\Aut(BS(p,q))$ is not finitely generated if $p$ properly divides $q$
(Theorem \ref{th:nfg}).

The key construction is an invariant deformation retract within
deformation space of $G$--trees.  This is the second purpose of this
paper.  The definition of a deformation space (of $G$--trees) appears
in the next section, but loosely speaking a deformation space \DS \ is
a moduli space of certain tree actions for a finitely generated group.
These spaces were introduced by Forester \cite{ar:F02}.  Culler and
Vogtmann's Outer space is a celebrated example of a deformation space
\cite{ar:CV86}.  Following intuition from Outer space, we define a
deformation retract $W \subset \DS$ of a general deformation space.
When the deformation space is Outer space, the deformation retract $W$
coincides with the spine of reduced Outer space.  In some cases we can
prove that $W$ is finite dimensional (Theorem \ref{th:fd}) and we
provide an example of a deformation space $\DS$ for which the
deformation retract $W$ is not finite dimensional (Example
\ref{ex:infinite-d}).  We note that Guirardel and Levitt have defined
a similar deformation retract within a non-ascending deformation space
\cite{ar:GL072}.

A non-solvable Baumslag--Solitar group $BS(p,q)$ has a canonical
deformation space $\DS_{p,q}$.  This deformation space is invariant
under the action of $\Out(BS(p,q))$.  The deformation retract of
$\DS_{p,q}$ mentioned in the preceding paragraph is denoted $W_{p,q}$.
We use this deformation space to construct the tree $X_{p,q}$.  If $p$
does not divide $q$ or $p = |q|$, then $W_{p,q}$ is a single point and
we take $X_{p,q} = W_{p,q}$.  If $|q|/p$ is prime, then the
deformation retract $W_{p,q}$ is a tree and we set $X_{p,q} =
W_{p,q}$.  In the other cases the complex $W_{p,q}$ is not a tree.
Using our description of $W_{p,q}$ when $|q|/p$ is prime, we define an
$\Out(BS(p,q))$--invariant subcomplex of $W_{p,q}$ and prove that this
subcomplex is a tree (Theorem \ref{th:Xtree}).  In this case we take
$X_{p,q}$ to be this tree.

We are hopeful that these techniques can extend to computing the
finiteness properties of other outer automorphism groups using
deformation spaces.

\begin{ack} 
  The research for this article was done under the supervision of my
  adviser Mladen Bestvina, whom I thank along with the University of
  Utah.  I am indebted to Vincent Guirardel and Gilbert Levitt for
  discussing this research with me and suggesting Theorem \ref{th:fd}.
  Also, I thank the referee for his/her helpful suggestions.
\end{ack}

%%%%%%%%%%%%%%%%%%%%%%%%%%%%%%%%%%%%%%%%%%%%%%%%%%%%%%%%%%%%%%%%%%%%%%%%%%%%% 
%\newpage %hard
\section{Deformation spaces of $G$--trees}\label{sc:ds}

%%%%%%%%%%%%%%%%%%%%%%%%%%%%%%%

\subsection{Definitions}\label{ssc:pre}

In this section we define the preliminary notions essential to the
following. 

For a graph $\Gamma$, we denote by $V(\Gamma)$ the set of vertices of
$\Gamma$ and by $E(\Gamma)$ the set of oriented edges of $\Gamma$.
For an edge $e \in E(\Gamma)$, $\bar{e}$ denotes the edge $e$ with the
opposite orientation.  The attaching maps are $\iv,\tv\co E(\Gamma)
\to V(\Gamma)$ (originating (initial) and terminating vertices).  For
a vertex $v \in V(\Gamma)$, let $E_\iv(v) = \{ e \in E(\Gamma) \ | \
\iv(e) = v\}$.  An \textit{edge path} $\gamma = (e_{0},\ldots,e_{n})$
is a set of oriented edges such that $\tv(e_{s}) = \iv(e_{s+1})$ for
$s = 0,\ldots,n-1$.  A \textit{circuit} is an edge path that is
homeomorphic to a circle and a \textit{loop} is a circuit consisting
of a single edge.

Let $G$ be a finitely generated group.  Later we will specialize to
Baumslag--Solitar groups, but the results and constructions in this
first section apply to a general finitely generated group.

A \textit{$G$--tree} $T$ is a simplicial tree that admits an action of
$G$ by simplicial automorphisms, i.e., by maps $f\co T \to T$ that are
bijections on the sets of edges and vertices and such that $f(\iv(e))
= \iv(f(e))$ for all $e \in E(T)$.  Actions are always assumed to be
without inversions, i.e., $ge \neq \bar{e}$ for all $g \in G, e \in
E(T)$.  Two $G$--trees are equivalent if there is a $G$--equivariant
simplicial isomorphism between them.  A \textit{metric $G$--tree} is a
$G$--tree with a metric such that the action of $G$ is by isometries.
As such we will consider two metric $G$--trees equivalent if there is
a $G$--equivariant isometry between them.  In either case, we will
always assume that the the $G$--tree $T$ is \textit{minimal} (no
invariant subtree), thus $T/G$ is a finite graph.

% It is impossible to talk about $G$--trees without using graphs of
% groups as the two notions are intimately related.  As such, the
% notation can become overwhelming and obtrusive.  We will try to keep
% things readable by maintaining some consistency.  In our examples,
% when we look at the set of edges emanating from a vertex $v$, edges in
% the same $G_{v}$--orbit will be denoted with the same letter, using a
% subscript to distinguish between them.  Several definitions for edges
% are made throughout the following.  These notions are always
% equivariant with respect to the $G$ action.  Therefore if some edge in
% the $G$--tree has a some property $\mathcal{P}$, then it also makes
% sense to say that the image of that edge in $T/G$ has property
% $\mathcal{P}$.  Some definitions are easier to state in terms of
% $T/G$.  In this case, an edge or edge path in $T$ has some property
% $\mathcal{P}$ if its image in $T/G$ does.

For a metric $G$--tree $T$, the \textit{length function}, $l_{T}\co G
\to \R$, is defined as $l_{T}(g) = \inf_{x \in T} d_{T}(x,gx)$.  It is
well-known that this infimum is achieved.  When we speak of a length
function for a simplicial $G$--tree $T$, we mean the length function
when we consider $T$ as a metric $G$--tree where all edges are
assigned length one.  The \textit{characteristic subtree} for an
element $g \in G$ is $T_{g} = \{ x \in T \ | \ d_{T}(x,gx) = l_{T}(g)
\}$.  If $l_{T}(g) > 0$, then $T_{g}$ is isometric to \R \ and $g$
acts on $T_{g}$ as a translation by $l_{T}(g)$.  In this case, $g$ is
called \textit{hyperbolic} and $T_{g}$ is called the \textit{axis} of
$g$.  Culler and Morgan proved that irreducible $G$--trees are
uniquely determined by their length functions \cite{ar:CM87}.  A
$G$--tree is \textit{irreducible} if $G$ does not fix an end of $T$
nor a pair of ends.  An equivalent condition is that there are two
hyperbolic elements whose axes are either disjoint or intersect in a
compact set \cite{ar:CM87}.

For a given $G$--tree $T$, a subgroup $H \subseteq G$ is
\textit{elliptic} if $H$ fixes a point in $T$.  There are two ways to
modify a $G$--tree, called \textit{collapse} and \textit{expansion},
that do not change the set of elliptic subgroups.  These two moves
correspond to the graph of groups isomorphism $A *_C C \cong A$.  A
finite sequence of these moves is called an \textit{elementary
deformation}.  Conversely, Forester proved that any two $G$--trees
with the same set of elliptic subgroups are related by an elementary
deformation \cite{ar:F02}.  The definitions of collapse and expansion
are as follows.

\begin{definition}[\textit{Collapse}]\label{def:co}
  Let $T$ be a $G$--tree and suppose there is a vertex $v \in V(T)$
  and an edge $e \in E_o(v)$ such that $G_e = G_v$ and $\iv(e)$ is not
  $G$--equivalent to $\tv(e)$, i.e., the image of $e$ in $T/G$ is not
  a loop.  We define a new $G$--tree $T_e$ by removing the edges $Ge$,
  then for all edges $f \in E(T)$ such that $\iv(f) = \tv(ge)$ for
  some $g \in G$ redefine $\iv(f) = gv$.  This is a \textit{collapse}
  move.  Such edges $e$ are called \textit{collapsible}.  
\end{definition}

If a $G$--tree does not admit a collapse move, it is called
\textit{reduced}.

\begin{definition}[\textit{Expansion}]\label{def:ex}
  Let $T$ be a $G$--tree and $v \in V(T)$.  Given a subgroup $H
  \subseteq G_v$ and an $H$--invariant set of edges $S \subseteq
  E_o(v)$ such that $G_e \subseteq H$ for all $e \in S$ we can define
  a new $G$--tree $T^{H,S}$ by adding a new edge $i$ via $\iv(i) = v$
  and redefining $\iv(e) = \tv(i)$ for $e \in S$, then repeating for
  all cosets $gH \in G/H$ using with the subgroup $gHg^{-1} \subseteq
  G_{gv}$ and the set of edges $gS \subseteq E_o(gv)$.  This is an
  \textit{expansion} move.
\end{definition}

The stabilizer of the new edge $i$ and the vertex $\tv(i)$ is $H$.

%%%%%%%%%%%%%%%%%%%%%%%%%%%%%%%%%%%%%%%%%%%%%%%%%%%%%%%%%%%%%%%%%%%%%%%%%%
%                                                                        %
%   \begin{figure}[h]                                                    %
%   \centering                                                           %
%   \input{exp_col.eepic}                                                %
%   \caption{Expansion and collapse in $G$--trees.}\label{fg:elem-def}   %
%   \end{figure}                                                         %
%                                                                        %
%%%%%%%%%%%%%%%%%%%%%%%%%%%%%%%%%%%%%%%%%%%%%%%%%%%%%%%%%%%%%%%%%%%%%%%%%%

% We give a more combinatorial description of an expansion move in
% Section \ref{ssc:ideal}.  An example of these moves in a $G$--tree is
% presented in Example \ref{ex:pre-i}.  

There are three special elementary deformations that we use
frequently: \textit{induction}, \textit{slide} and $\As^{\pm 1}$.
Pictures for these moves in the general setting appear in
\cite{ar:CFpp}.  In Section \ref{sc:gbs} these moves are shown for the
canonical deformation space associated to a generalized
Baumslag--Solitar group.

\begin{definition}[\textit{Induction}]\label{def:in}
  Let $T$ be a $G$--tree and suppose that at the vertex $v \in V(T)$
  there are two edges $e,e' \in E_o(v)$ such that $G_e = G_v$, $G_{e'}
  \neq G_v$ and $e' \in G\bar{e}$, i.e., $e$ and $e'$ project to the
  same loop in $T/G$ with opposite orientations.  The composition of
  the expansion using any subgroup $G_{e'} \subseteq H \subsetneq G_v$
  and the set $S = He'$ followed by the collapse of $e$ is called an
  \textit{induction}.  The inverse of an induction is also called an
  induction.
\end{definition}

The loop in $T/G$ created by the image of $e$ in the above definition
is called an \textit{ascending loop}, i.e.~a loop in $T/G$ in which
one of the attaching maps is an isomorphism and the other is a proper
inclusion.  An induction move changes the stabilizer of the vertex $v$
to $H$.

% \begin{figure}[h]
% \centering
% \input{induction.eepic}
% \caption{Induction in a $G$--tree.}\label{fg:induct}
% \end{figure}

\begin{definition}[\textit{Slide}]\label{def:sl}
  Let $T$ be a $G$--tree and suppose that at the vertex $v \in V(T)$
  there are two edges $e,f \in E_o(v)$ such that $G_f \subseteq G_e$
  and $f$ is not $G$--equivalent to $e$ or $\bar{e}$.  We can perform
  an expansion using the subgroup $H = G_e$ and set $S = \{e\} \cup
  G_ef$.  Then in the expanded $G$--tree, the edge $e$ is collapsible.
  The composition of this expansion followed by the collapse of $e$ is
  called a \textit{slide}.
\end{definition}  

% \begin{figure}[h]      
% \centering             
% \input{slide.eepic}
% \caption{Sliding the edge $f$ over $e$.\label{fg:slide}}   
% \end{figure}

A slide can also be thought of as removing the edge $f$ at $\iv(f)$
and reattaching it via $\iv(f) = \tv(e)$, then repeating equivariantly
throughout $T$.  Topologically, we are folding the first half of $f$
across $e$.

\begin{definition}[\textit{$\As^{\pm 1}$--move}]\label{def:A}
  Let $T$ be a $G$--tree and suppose that at the vertex $v \in V(T)$
  there are edges $e,e',f \in E_o(v)$ where $e,e'$ satisfy the
  hypotheses of an induction move and $G_{e'} \subseteq G_f$.  Further
  suppose that there are exactly three $G_v$--orbits in $E_o(v)$ and
  the edges $e,e'$ and $f$ are in distinct orbits.  Then there is an
  induction move after which $G_v = G_f$.  The composition of this
  induction move followed by the collapse of the edge $f$ is called an
  $\As^{-1}$--move.  The inverse is called an $\As$--move.
\end{definition}

An $\As$--move creates an ascending loop and an $\As^{-1}$--move
removes an ascending loop.  
% As in \cite{ar:CFpp}, we require that both
% the initial $G$--tree and the $G$--tree resulting from either of the
% abve three moves to be reduced.

Deformation spaces of $G$--trees were introduced by Forester
\cite{ar:F02}.  Given a $G$--tree $T$, let $\uDS$ be the set of all
metric $G$--trees that define the same set of elliptic subgroups as
$T$.  This set of metric $G$-trees is called an \textit{unnormalized
  deformation space of $G$--trees}.  There is an action of $\R^+$ on
$\uDS$ by scaling the metric on a given metric $G$--tree.  The
quotient $\DS$, is called a \textit{deformation space of $G$--trees},
or sometimes just a \textit{deformation space}.  By Forester's
deformation theorem, disregarding the metric on any two projectivized
$G$--trees $T,T' \in \DS$, there is an elementary deformation
transforming $T$ to $T'$.

The importance of the three special moves described above is given
by the following theorem.

\begin{theorem}[\cite{ar:CFpp}]\label{th:ais}
  In a deformation space of $G$--trees, any two reduced trees are
  related by a finite sequence of slides, inductions and $\As^{\pm
    1}$--moves, with all intermediate trees reduced.
\end{theorem}

Deformation spaces can be topologized in several ways: \textit{axes
  topology, equivariant Gromov--Hausdorff topology} or \textit{weak
  topology}.  In the case of irreducible $G$--trees, the axes topology
and the equivariant Gromov--Hausdorff topology are the same
\cite{ar:P89}.  In addition for locally finite $G$--trees, the
equivariant Gromov--Hausdorff topology and the weak topology are the
same \cite{ar:GL072}.  This is not true in general however, for an
example see \cite{ar:MM96}.  The topology we work with is the weak
topology; for definitions of the other two see \cite{ar:C05}.  A
projectivized $G$--tree $T \in \DS$ determines an open simplex by
equivariantly changing the lengths of the edges of $T$ while holding
the sum of the lengths of the edges in $T/G$ constant.  This open
simplex has dimension one less than the number of edges of $T/G$.  The
faces in the closure of this open simplex are found by collapsing
subsets of collapsible edges in $T$.

Culler and Vogtmann's Outer space is an example of a deformation space
of $G$--trees \cite{ar:CV86}.  In this case, $G$ is a finitely
generated free group of rank $n$ at least 2, $F_{n}$, and the only
elliptic subgroup is the trivial group.  In other words, all of the
actions are free.  Another example is the complex $K_0(G)$ defined by
McCullough and Miller \cite{ar:MM96}.  Here, $G$ is an arbitrary
finitely generated group and the elliptic subgroups are the free
product factors in a Grusko decomposition for $G$.  Guirardel and
Levitt describe a different complex for a free product decomposition
of a finitely generated group $G$ \cite{ar:GL07}.  Their complex is
the deformation space where the elliptic subgroups are the free
product factors in a Grusko decomposition for $G$ that are not
infinite cyclic.

The deformation space $\DS$ is acted upon by some subgroup
$\Out(G)_{\DS} \subseteq \Out(G)$.  This is the subgroup that
preserves the set of conjugacy classes of the elliptic subgroups
associated to $\DS$. In the case of Culler and Vogtmann's Outer space,
this is the entire group $\Out(F_{n})$.  Likewise, the complex defined
by Guirardel and Levitt is invariant under $\Out(G)$.  For the
McCullough--Miller complex $K_0(G)$, $\Out_{\DS}(G) = \SOut(G)$, the
subgroup of \textit{symmetric outer automorphisms}.  These examples
are contractible and have been used for computing some of the
finiteness properties of $\Out_{\DS}(G)$, see individual references.

Under some mild hypotheses, we \cite{ar:C05} and independently
Guirardel and Levitt \cite{ar:GL072} have shown that deformation
spaces are contractible.  Both our proof and the proof of Guirardel
and Levitt use Skora's method of continuous folding \cite{ar:Spp}.

\begin{theorem}\label{th:cont}
For a finitely generated group $G$, any irreducible deformation space
that contains a $G$--tree with finitely generated vertex stabilizers is
contractible.
\end{theorem}

A deformation space is \textit{irreducible} if all $G$--trees
(equivalently a single $G$--tree) in $\DS$ are (is) irreducible.  In
\cite{ar:C05}, the above theorem is only shown for the equivariant
Gromov--Hausdorff topology (equivalently axes topology).  Guirardel
and Levitt show this as well as showing that the contraction is
continuous in the weak topology.  In addition, they replace
irreducible with a weaker hypothesis.

%%%%%%%%%%%%%%%%%%%%%%%%%%%%%%%

\subsection{A deformation retract of \DS}\label{ssc:def}

Culler and Vogtmann's Outer space deformation retracts to a subcomplex
called \textit{reduced Outer space}.  This is the subcomplex of
projectivized free metric $F_{n}$--trees $T$ such that the quotient
graph $T/F_{n}$ does not have any separating edges.  This deformation
retract is obtained by equivariantly shrinking the edges that project
to separating edges in the quotient graphs $T/F_{n}$.  There is a
further deformation retraction of reduced Outer space to a spine.  We
will describe a similar deformation retract of a deformation space
$\DS$.  For general deformation spaces, it is easier to deformation
retract to the spine first, then define an additional deformation
retraction.

Let $K$ be the spine of $\DS$.  Specifically, let $\mathcal{OS}(\DS)$
be the poset of open simplices of $\DS$ where $\sigma \leq \sigma'$ if
$\sigma$ is a face in the closure of $\sigma'$.  The spine $K$ is
defined as the geometric realization of the poset $\mathcal{OS}(\DS)$.
The spine $K$ is an $\Out(G)_{\DS}$--invariant deformation retraction
of the deformation space $\DS$ (with the weak topology).  Therefore,
$\Out(G)_{\DS}$ acts on $K$, which is a contractible simplicial
complex.

The spine $K$ has an alternative combinatorial description in terms of
the $G$--trees appearing in $\DS$. The poset $\mathcal{OS}(\DS)$ is
isomorphic to the poset $\mathcal{C}ol(\DS)$ of $G$--trees in $\DS$
(thought of only as simplicial trees) where $T' \geq T$ if $T'$
collapses to $T$.  We say that $T'$ \textit{collapses to} $T$ if there
is a $G$--equivariant map $T' \to T$ that is a composition of finitely
many collapse moves.  Thus $K$ can also be viewed as the geometric
realization of the poset $\mathcal{C}ol(\DS)$.  Vertices of $K$ are
$G$--trees in $\DS$; higher dimensional simplices correspond to
collapse sequences of such $G$--trees.  To ease notation, we will use
$K$ to denote the set of vertices of $K$.  For future reference, we
remark that no $G$--tree $T \in K$ contains a \textit{subdivision
  vertex}, i.e.~a vertex $v \in V(T)$ with $E_o(v) = \{e,f\}$ and $G_e
= G_v = G_f$.

For a $G$--tree $T \in K$ define the following sets:
\begin{align*}
  \col(T) &= \{ T' \in K \ | \ T \geq T' \}, \\
  \red(T) &= \{T' \in \col(T) \ | \ T' \mbox{ is reduced} \}, \mbox{
    and} \\
  \mathcal{W}(T) &= \{ T' \in \col(T) \ | \ \red(T) = \red(T') \}.
\end{align*}
In the setting of Culler and Vogtmann's Outer space, $\red(T)$ is the
set of $F_{n}$--trees found by collapsing some maximal forest of
$T/F_{n}$ and $\mathcal{W}(T)$ is the set of $F_{n}$--trees found by
collapsing separating edges in $T/F_{n}$.  The following lemma shows
that the above sets can be realized as certain subsets of collapsible
edges in $T/G$.

\begin{lemma}\label{lm:neq}
  If $T$ is an irreducible $G$--tree with collapsible edges $e$ and
  $f$ where $f \notin Ge \cup G\bar{e}$, then $T_{e} \neq T_{f}$.
  Also, if $e$ is collapsible in $T_{f}$ and $f$ is collapsible in
  $T_{e}$, then $(T_{e})_{f} = (T_{f})_{e}$.
\end{lemma}

\begin{proof}
  If there is a hyperbolic element $g \in G$ whose axis projects down
  to a closed path in $T/G$ that crosses the image of $e$ more than it
  crosses the image of $f$ then $l_{T_{e}}(g) < l_{T_{f}}(g)$.  Thus
  as irreducible $G$--trees are determined by their length functions,
  $T_{e} \neq T_{f}$.  It is easy to find such an element $g \in G$ by
  looking at edge paths (in the graph of groups sense) in $T/G$, see
  \cite{ar:B93,ar:S80}.

  By looking at length functions again, the second part of the lemma
  is obvious.
\end{proof}

The following definition appears in \cite{ar:GL072}.

\begin{definition}\label{def:sur}
  An edge $e \in E(T)$ is called \textit{surviving} it there is a
  $G$--tree $T' \in \red(T)$ such that $e$ is not collapsed in $T \to
  T'$.  If an edge is not a surviving edge, it is called
  \textit{non-surviving}.
\end{definition}

\begin{lemma}\label{lm:T_{W}}
  $\mathcal{W}(T)$ has a unique minimal element, $T_{\mathcal{W}}$.
  This $G$--tree is characterized as the $G$--tree obtained from $T$
  by collapsing the non-surviving edges.  Further, if $T \geq T'$ then
  $T_{\mathcal{W}} \geq T'_{\mathcal{W}}$.
\end{lemma}    
    
\begin{proof}
  As $\mathcal{W}(T)$ is a finite set, minimal elements exist in
  $\mathcal{W}(T)$.  Let $T_{0}$ be a minimal element.  Then every
  non-surviving edge in $T$ must be collapsed in $T \to T'_{0}$ since
  $T_{0}$ is minimal.  Also, no surviving edge could be collapsed in
  $T \to T_{0}$ as $\red(T) = \red(T_{0})$.  Hence any minimal element
  $T_{0}$ is found by collapsing the non-surviving edges.  By Lemma
  \ref{lm:neq}, this completely determines $T_{0}$.
    
  Finally notice that if $T \geq T'$ then any non-surviving edge in
  $T$ is either collapsed in $T \to T'$ or it is non-surviving for
  $T'$.  This is true since if the edge is surviving for $T'$, then it
  must also be surviving for $T$.  Thus $T_{\mathcal{W}} \geq
  T'_{\mathcal{W}}$.
\end{proof}

Define $h\co \mathcal{C}ol(\DS) \to \mathcal{C}ol(\DS)$ by $h(T) =
T_\mathcal{W}$.  The above lemma shows that $h$ is a well-defined
poset map.  Notice that $h(T) \leq T$.  The following lemma shows that
$h$ defines a deformation retract of $K$, the geometric realization of
$\mathcal{C}ol(\DS)$.

\begin{lemma}[Quillen's Poset Lemma \cite{ar:Q78}]\label{lm:poset}
  Let $X$ be a poset and $f\co X \to X$ be a poset map (i.e. $x \leq
  x'$ implies $f(x) \leq f(x')$ for all $x,x' \in X$) with the
  property that $f(x) \leq x$ for all $x \in X$ (or $f(x) \geq x$ for
  all $x \in X$).  Then the geometric realization of $f(X)$ is a
  deformation retract of the geometric realization of $X$.
\end{lemma}

Define $W$ as the geometric realization of the poset
$h(\mathcal{C}ol(\DS))$.  Therefore, by Quillen's Poset Lemma, $W$ is
a $\Out(G)_{\DS}$--invariant deformation retract of \DS. In
particular, $W$ is contractible. For Culler and Vogtmann's Outer
space, $W$ is the spine of reduced outer space.  We will explicitly
show the deformation retraction $h\co K \to W$ in Example
\ref{ex:star} for a star in a deformation space for the
Baumslag--Solitar group $BS(2,4)$.  If $\DS$ is a non-ascending
deformation space, then $W$ is the spine of the deformation retract
defined by Guirardel and Levitt \cite{ar:GL072}.  A deformation space
$\DS$ is called \textit{non-ascending} if for all $ T \in \DS$, the
quotient graph of groups $T/G$ does not contain an ascending loop.  As
for $K$, we will use $W$ to denote the set of vertices of $W$.

\subsection{$G$--trees in $W$}\label{ssc:lkW}

In this section we determine when a $G$--tree in $\DS$ represents a
vertex in $W$.  The following definition generalizes the definition of
\textit{shelter} in \cite{ar:GL072}.

\begin{figure}[ht]
\centering
\psfrag{S1}{S1}
\psfrag{S2}{S2}
\psfrag{S3}{S3}
\psfrag{neq}{{\small$\neq$}}
\psfrag{=}{{\small$=$}}
\psfrag{cdots}{$\cdots$}
\psfrag{o}{$\cdot$}
\includegraphics{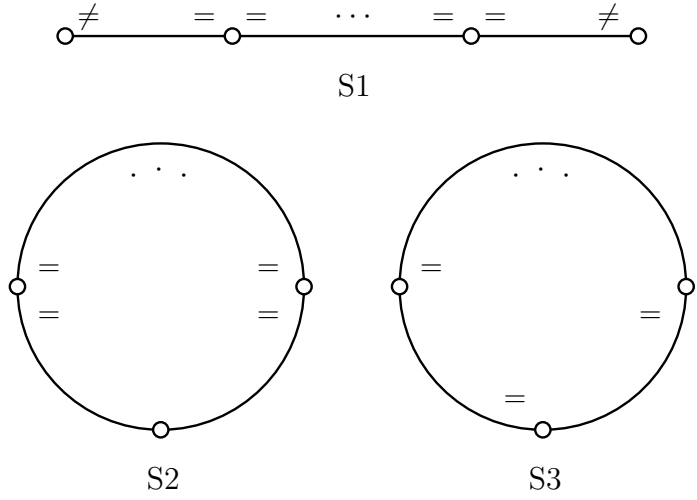} 
\caption{Shelters in $T/G$.  The ``='' signs denote when an inclusion
  $G_{e} \subseteq G_{v}$ is an equality.  The ``$\neq$'' signs denote
  when an inclusion is proper.  Other inclusions may or may not be
  proper.}\label{fg:shelters}
\end{figure}

\begin{definition}\label{def:shelter}
  Let $T$ be a $G$--tree, $\gamma = (e_{0},\ldots,e_{n}) \subseteq T$
  an edge path in $T$ and $\widehat{\gamma}$ the image of $\gamma$ in
  $T/G$.  We say $\gamma$ (or $\widehat{\gamma}$) is a
  \textit{shelter} if either:
\begin{itemize}
    
\item[(S1)] $\widehat{\gamma}$ is a topological segment,
  $G_{\iv(e_{0})} \neq G_{e_{0}}$, $G_{e_{n}} \neq G_{\tv(e_{n)}}$ and
  $G_{e_{s}} = G_{\tv(e_{s})} = G_{e_{s+1}}$ for $s = 0, \ldots, n-1$;

\item[(S2)] $\widehat{\gamma}$ is a circuit and $G_{e_{s}} =
  G_{\tv(e_{s})} = G_{e_{s+1}}$ for $s = 0,\ldots,n-1$; 

\item[(S3)] $\widehat{\gamma}$ is a circuit and $G_{\iv(e_{s})} =
  G_{e_{s}}$ for $s = 0,\ldots,n$.

\end{itemize}
We refer to the labels S1, S2 or S3 as the \textit{type} of the
shelter.  See Figure \ref{fg:shelters}.
\end{definition}

The following proposition generalizes Corollary 7.5(1) in \cite{ar:GL072}.

\begin{proposition}\label{prop:shelter}
  An edge is surviving if and only if it is contained in a shelter.
\end{proposition}

Before we prove this, we prove a simple lemma about how stabilizers
can change after collapse moves.

\begin{lemma}\label{lm:stabilizer}
  Let $f$ be a collapsible edge in $T$.  Let $\iv(f) = v, \tv(f) = w$
  and denote the image of $f$ in $T_{f}$ by $z$.  Then for $e \in
  E(T)$, $e \notin Gf \cup G\bar{f}$ with $\iv(e) = v$, $G_{e} \neq
  G_{z}$ if and only if $G_{e} \neq G_{v}$ or $G_{f} \neq G_{w}$.
  %, in which case $G_{f} = G_{v}$.
\end{lemma}
    
\begin{proof}
  As $f$ is collapsible, $G_{v} \subseteq G_{w}$ or $G_{w} \subseteq
  G_{v}$, in either case $G_{z}$ is the union of the two subgroups.
  Therefore if $G_{e} \subsetneq G_{v}$, then $G_{e} \neq G_{z}$.  If
  $G_{f} \subsetneq G_{w}$, then as $f$ is collapsible $G_{f} =
  G_{v}$.  Therefore $G_{e} \subseteq G_{v} = G_{f} \subsetneq G_{w} =
  G_{z}$.  In particular, $G_{e} \neq G_{z}$.

  For the converse, suppose that $G_{e} = G_{v}$ and $G_{f} = G_{w}$.
  Then clearly $G_{z} = G_{v}$, hence $G_{e} = G_{z}$.
\end{proof}

Now we can prove Proposition \ref{prop:shelter}.

\begin{proof}[Proof of Proposition \ref{prop:shelter}]
  It is obvious that every edge in a shelter is surviving.  The
  converse is straight forward, but there are several cases.  Some of
  these cases are presented in \cite{ar:GL072}, but we present an
  entire proof for completeness.  Let $T$ by a $G$--tree and suppose
  that $e$ is a surviving edge in $T$.  Let $T'$ be a reduced
  $G$--tree such that $e$ is not collapsed in $T \to T'$.  We denote
  the image of $e$ in $T/G$ by $e$.

  \medskip \noindent {\bf case 1}: The image of $e$ in $T'/G$ is
  an interval.

  Let $Y$ be the maximal subtree of $T/G$ the contains $e$ and is
  collapsed to $e$ in $T/G \to T'/G$.  We will show by induction on
  the number of edges in $Y$ that $Y$ contains a shelter of type S1
  that contains $e$.  If the number of edges in $Y$ is one, then $Y=e$
  and as $T'$ is reduced, $e$ is a shelter.

  Now suppose that the number of edges in $Y$ is greater than one.
  Let $f$ be an edge of $Y$ other than $e$.  Then the image of $Y$ in
  $T_{f}/G$ is the maximal subtree in $T_{f}/G$ that contains $e$ and
  collapses to $e$ in $T_{f}/G \to T'/G$.  Hence, by induction the
  image of $Y$ in $T_{f}/G$ contains a shelter $(e_{0},\ldots,e_{n})$
  of type S1 containing $e$.  We consider these edges as edges of
  $T/G$.

  Suppose that $f$ is adjacent to $\iv(e_{0})$.  Orient $f$ such that
  $\tv(f) = \iv(e_0)$.  Then by Lemma \ref{lm:stabilizer}, either
  $(e_{0},\ldots,e_{n})$ or $(f,e_{0},\ldots,e_{n})$ is a shelter of
  type S1 for $e$ as $(e_{0},\ldots,e_{n})$ is a shelter of type S1 in
  $T_{f}/G$.  Similarly, we can find a shelter containing $e$ of type
  S1 if $f$ is adjacent to $\tv(e_{n})$.

  Next suppose that $f$ is adjacent to $\tv(e_{s})$ and $\iv(e_{s+1})$
  for $0 \leq s \leq n-1$.  If $\tv(e_{s}) = \iv(e_{s+1})$ in $T/G$,
  then by Lemma \ref{lm:stabilizer} $(e_{0},\ldots,e_{n})$ is a
  shelter of type S1 for $e$.  If not, then orient $f$ such that
  $\tv(e_{s}) = \iv(f)$ and $\iv(e_{s+1}) = \tv(f)$.  Then again by
  Lemma \ref{lm:stabilizer}
  $(e_{0},\ldots,e_{s},f,e_{s+1},\ldots,e_{n})$ is a shelter of type
  S1 for $e$.

  \medskip \noindent {\bf case 2}: The image of $e$ in $T'/G$ is a
  loop.

  Let $Y$ be the circuit in $T/G$ that collapses to $e$ in $T/G \to
  T'/G$.  We again use induction on the number of edges in $Y$.  Our
  claim in this case is that either $Y$ is a shelter of type S2 or S3
  or else $Y$ contains part of a shelter of type S1 that contains $e$.
  If $Y$ contains only one edge, then $Y=e$ and $Y$ is a shelter of
  type S2 or S3 as $T'/G$ is reduced.  Otherwise as before, take any
  other edge $f$ in $Y$ other than $e$ and look at the image of the
  circuit $Y$ in $T_{f}/G$.  This is the circuit in $T_f/G$ that
  collapses to $e$ in $T'/G$.

  If the image of $Y$ in $T_{f}/G$ contains part of a shelter of type
  S1 that contains $e$, then proceed as in case 1 to show that $Y$ in
  $T/G$ contains a shelter of type S1 that contains $e$.

  Otherwise, we suppose that the image of $Y$ in $T_{f}/G$ is a
  shelter of type S2 or S3.  Suppose that $\iv(f) = \tv(e_{0}) =
  \iv(e_{1})$ for two edges $e_{0},e_{1} \subseteq Y$.  Then it is
  simple check using Lemma \ref{lm:stabilizer} that $Y$ is shelter of
  type S2 or S3.

  Finally, suppose that in $T/G$ we have $\tv(e_{0}) = \iv(f)$ and
  $\iv(e_{1}) = \tv(f)$.  Denote the image of $f$ in $T_{f}/G$ by $z$.
  First we assume that the image of $Y$ is a shelter of type S2.  If
  $G_{e_{0}} = G_{z} = G_{e_{1}}$ then by Lemma \ref{lm:stabilizer},
  $G_{e_{0}} = G_{\iv(f)} = G_{f}$ and $G_{\iv(e_{1})} = G_{\tv(f)} =
  G_{f}$, hence $Y$ is a shelter of type S2.  If $G_{e_{0}} \neq G_{z}
  = G_{e_{1}}$ then either $G_{e_{0}} \neq G_{\iv(f)}$ in which case
  $Y$ is a shelter of type S2 or S3 (depending on whether $G_{f} =
  G_{\tv(f)}$ or $G_{f} = G_{\iv(f)}$) or $G_{e_{0}} = G_{\iv(f)}$ and
  $G_{f} \neq G_{\tv(f)}$, in which case $Y$ is shelter of type S2.
  If $G_{e_{0}} \neq G_{z} \neq G_{e_{1}}$ then if $G_{e_{0}} \neq
  G_{\iv(f)}$ and $G_{e_{1}} \neq G_{\tv(f)}$ then $Y - \{f\}$ is a
  shelter of type S1 that contains $e$.  Otherwise suppose that
  $G_{e_{1}} = G_{\tv(f)}$, then by a similar argument as before, $Y$
  is a shelter of type S2.  Now assume that the image of $Y$ in
  $T_{f}/G$ is a shelter of type S3.  If $G_{e_{0}} = G_{z} =
  G_{e_{1}}$, then by Lemma \ref{lm:stabilizer}, $Y$ is a shelter of
  type S3.  If $G_{e_{0}} = G_{z} \neq G_{e_{1}}$, then $G_{\tv(f)} =
  G_{f}$ and $Y$ is a shelter of type S3.

  This completes the proof of Proposition \ref{prop:shelter}.
\end{proof}

As a corollary, we get a condition to check whether or not a $G$--tree
is in $W$.

\begin{corollary}\label{co:tw=t}
  $T = T_{\mathcal{W}}$ and hence $T \in W$ if and only if $T/G$
  is a union of shelters.
\end{corollary}

\begin{proof}
  This follows immediately from Lemma \ref{lm:T_{W}} and Proposition
  \ref{prop:shelter} as $T = T_{\mathcal{W}}$ if and only if every
  edge in $T$ is surviving.
\end{proof}

% Later on we will give an example of $\str_{K}(T),\str_{W}(T)$ for a
% $BS(2,4)$--tree (Example \ref{ex:star}).  Before we do so, we will
% give two sufficient conditions guaranteeing that the complex $W$ is
% finite dimensional.

%%%%%%%%%%%%%%%%%%%%%%%%%%%%%%%

\subsection{Finite Dimensionality of $W$}\label{ssc:fdW}

We conclude our treatment of deformation spaces for a general finitely
generated group with a discussion of the finite dimensionality of the
deformation retract $W$.

If $T$ and $T'$ are $G$--trees in $W$ then there is an elementary
deformation taking $T \to T'$.  As each of the elementary moves is a
homotopy equivalence of the quotient graph, $T/G$ is homotopy
equivalent to $T'/G$.  Therefore homotopy invariants of graphs are
invariants of deformation spaces. 

\begin{lemma}\label{lm:nafd}
  Let $\DS$ be a non-ascending deformation space for a finitely
  generated group $G$.  Then the number of vertices in $T/G$ for any
  $G$--tree $T \in W$ is bounded.
\end{lemma}

\begin{proof}
  By the above remark, the Euler characteristic $\chi(T/G)$ is
  constant for $T \in \DS$.  We can compute the Euler characteristic
  by $\chi(T/G) = \frac{1}{2}(V_{1} - V_{3} - 2V_{4} - 3V_{5} -
  \ldots)$, where $V_{s}$ denotes the number of vertices with valence
  $s$.  Therefore, as in \cite{ar:BF91}, it suffices to show that
  $V_{1}$ and $V_{2}$ are bounded.  Let $N$ denote the number of
  conjugacy classes of maximal elliptic subgroups of $G$.  This number
  is finite and depends only on $\DS$.

  By minimality of $T \in W$, every valence one vertex in $T/G$
  corresponds to a unique conjugacy class of a maximal elliptic
  subgroups.

  By Corollary \ref{co:tw=t}, $T/G$ is a union of shelters.  As $\DS$
  is non-ascending, only shelters of type S1 and S2 appear.  Further
  as $\DS$ is non-ascending if $e_0,e_1$ are adjacent edges in a
  shelter of type S2 with $\tv(e_0) = \iv(e_1) = v$ and $G_{e_0} \neq
  G_v$, then $G_{e_1} \neq G_v$ also.

  Let $v$ be a valence two vertex in $T/G$ with adjacent edges $e,f$.
  Suppose that $v$ is contained in a shelter of type S2.  Then since
  there are no subdivision vertices in $G$--trees in $\DS$, $G_e \neq
  G_v$ and $G_f \neq G_v$ by the above remark, hence $v$ corresponds
  to a unique conjugacy class of maximal elliptic subgroups.  If $v$
  is contained in a shelter of type S1, then as there are no
  subdivision vertices in $G$--trees in $W$, $v$ must be the endpoint
  of two shelters of type S1.  Hence again, $G_e \neq G_v$ and $G_f
  \neq G_v$ and $v$ corresponds to a unique conjugacy class of maximal
  elliptic subgroups

  This shows that $V_{1}+V_{2} \leq N$.  Therefore, the number of
  vertices in $T/G$ for any $G$--tree $T \in W$ is bounded.
\end{proof}

% We remark that the deformation retract of $\DS$ defined by Guirardel
% and Levitt when $\DS $ is non-ascending is finite dimensional
% \cite{ar:GL072}.
We will see that this above lemma implies that the deformation retract
$W \subset \DS$ for a non-ascending deformation space is finite
dimensional.  Before doing so, we look at a different setting where we
can bound the number of vertices in $T/G$.  As remarked above, the
first Betti number $b_1(T/G)$ of a quotient graph of $T \in \DS$
defines an invariant of the deformation space $\DS$, which abusing
notation, we denote by $b_1(\DS)$.  If $b_1(\DS) = 0$ then $\DS$ is
non-ascending as there are no loops in $T/G$ for any $T \in \DS$.

A deformation space is \textit{locally finite} if every $G$--tree $T
\in \DS$ is locally finite.  As modifying a locally finite $G$--tree
by an elementary deformation results in a locally finite tree, a
deformation space is locally finite if a single $G$--tree $T \in \DS$
is locally finite.  We remark that if $\DS$ is locally finite (as a
deformation space) then $\DS$ and the deformation retract $W$ are
locally finite as simplicial complexes although the converse is not
true.

Bass and Kulkarni introduced an invariant of a locally finite
deformation space $\DS$, called the \textit{modular homomorphism}
$q_\DS\co G \to \Q^\times$ \cite{ar:BK90}.  This homomorphism is
defined by:
\begin{equation*}
q_\DS(g) = [V : V \cap V^g] / [V^g : V \cap V^g]
\end{equation*}
where $V$ is any subgroup of $G$ commensurable to a vertex stabilizer
for a $G$--tree $T \in \DS$.  There is a useful alternative
description of this homomorphism.  For an edge $e \in E(T)$, define
$i(e) = [G_{\iv(e)} : G_e]$ and $q(e) = i(e) / i(\bar{e}) \in
\Q^\times$.  This map $q$ descends to edges in $T/G$ and hence also to
$H_1(T/G)$ by multiplication.  The homomorphism $q_\DS\co G \to
\Q^\times$ is the composition $G \to H_1(T/G) \to \Q^\times$
\cite{ar:F06}.

If $q_\DS(G) \cap \Z = \{1\}$ then $\DS$ is non-ascending as the
modulus of any ascending loop is a non-trivial integer.  The converse
is not true.  Forester showed that if $q_\DS(G) \cap \Z = \{1\}$ then
the canonical deformation space for a generalized Baumslag--Solitar
group $\DS$ is finite dimensional \cite{ar:F06}.  In fact Forester
showed that the quotient $W / \Out(G)_{\DS}$ is compact in this case.
We remark that the \Z--rank of the subgroup $q_\DS(G) \subseteq
\Q^\times$ is bounded by $b_1(\DS)$.

\begin{lemma}\label{lm:fdb1}
  Let $\DS$ be a locally finite irreducible deformation space for a
  finitely generated group with $b_1(\DS) = 1$.  Then the number of
  vertices in $T/G$ for any $G$--tree $T \in W$ is bounded.
\end{lemma}

\begin{proof}
  The proof is similar to the proof of Lemma \ref{lm:nafd}, we must
  show that there is a bound on the number of valence one and two
  vertices.  Again, let $N$ be the number of conjugacy classes of
  maximal elliptic subgroups.  As before, each valence one vertex
  corresponds to a unique conjugacy class of maximal elliptic
  subgroups.

  If all of the shelters are type S1, then every valence two vertex
  must be the endpoint of the two shelters it is in and hence
  corresponds to a unique conjugacy class of maximal elliptic
  subgroups.  Therefore, suppose that we have a shelter of type S2 or
  S3.  As $b_{1}(\DS) = 1$, there is only one such shelter.  Any
  valence two vertex not in this shelter must by in a shelter of type
  S1 and as above it corresponds to a unique conjugacy class of
  maximal elliptic subgroups.  Therefore, we only need to bound the
  number of valence two vertices in the shelter.

  If the shelter is type S2, then there can be at most one valence two
  vertex as there are no subdivision vertices in $G$--trees in $W$.
  Otherwise, if the shelter is type S3, then $q_{\DS}(G)$ is generated
  by an integer $q$ and the number of valence two vertices is bounded
  by the number of prime factors in $q$ as this bounds the length of a
  chain of proper subgroup inclusions that can appear in $T/G$.
\end{proof}

\begin{theorem}\label{th:fd}
  Let $G$ be a finitely generated group and $\DS$ an irreducible
  deformation space for $G$.  If $\DS$ is either:
\begin{itemize}
    
\item[1.] non-ascending; or
    
\item[2.] locally finite and has $b_1(\DS) \leq 1$;

\end{itemize}
then the deformation retract $W \subset \DS$ is finite dimensional.
\end{theorem}

\begin{proof}
  Recall that a simplex in $W$ is a sequence of collapse moves between
  $G$--trees in $W$.  As a collapse move reduces the number of
  vertices in the quotient graph by at least one, finite
  dimensionality of $W$ is equivalent to a uniform bound on the number
  of vertices in $T/G$ for any $T \in W$.  This is the content of
  Lemmas \ref{lm:nafd} and \ref{lm:fdb1}.
\end{proof}

The first part of the previous theorem also appears as Theorem 7.6 in
\cite{ar:GL072}.  An example of a deformation space $\DS$ for which
$W$ is not finite dimensional is presented in the next section
(Example \ref{ex:infinite-d}).

%%%%%%%%%%%%%%%%%%%%%%%%%%%%%%%%%%%%%%%%%%%%%%%%%%%%%%%%%%%%%%%%%%%%%%%%

\section{Deformation spaces for GBS groups}\label{sc:gbs}

A group $G$ that acts on a tree such that the stabilizer of any point
is infinite cyclic is called a \textit{generalized Baumslag--Solitar
  (GBS) group}.  These groups have also recently appeared in
\cite{ar:F06,ar:L07}.

\begin{figure}[h]      
\centering             
\psfrag{p}{{\small$p$}}
\psfrag{q}{{\small$q$}}
\psfrag{2}{{\small 2}}
\psfrag{3}{{\small 3}}
\psfrag{5}{{\small 5}}
\psfrag{4}{{\small 4}}
\psfrag{7}{{\small 7}}
\includegraphics{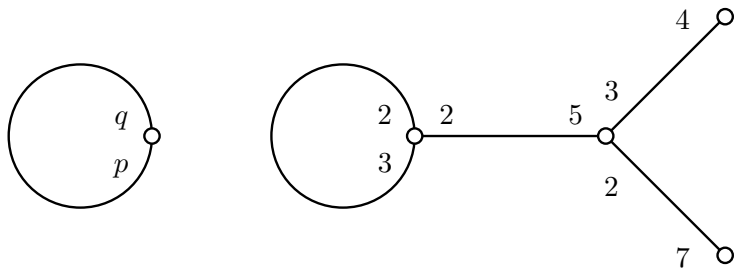} 
\caption{Examples of labeled graphs.  The labeled graph on the left
  represents the classical Baumslag--Solitar group
  $BS(p,q)$.}\label{fg:label}
\end{figure}

This action of a GBS group determines a graph of groups decomposition
of $G$ where all of the edge groups and vertex groups are isomorphic
to \Z. As such, all of the attaching maps are given by multiplication
by some nonzero integer.  This data can be represented succinctly as a
\textit{labeled graph}.  Specifically, a labeled graph is a pair
$(\Gamma,\lambda)$ where $\Gamma$ is a finite graph and $\lambda\co
E(\Gamma) \to \Z - \{0\}$ is a function.  The labels $\lambda(e)$
represent, for a chosen set of generators for the vertex and edge
groups, the inclusion maps $G_{e} \to G_{v}$.  See Figure
\ref{fg:label} for examples of labeled graphs.

\begin{figure}[ht]
\centering
\psfrag{expansion}{expansion}
\psfrag{collapse}{collapse}
\psfrag{a}{{\small$a$}}
\psfrag{b}{{\small$b$}}
\psfrag{n}{{\small$n$}}
\psfrag{1}{{\small 1}}
\psfrag{nc}{{\small$nc$}}
\psfrag{nd}{{\small$nd$}}
\psfrag{c}{{\small$c$}}
\psfrag{d}{{\small$d$}}
\includegraphics{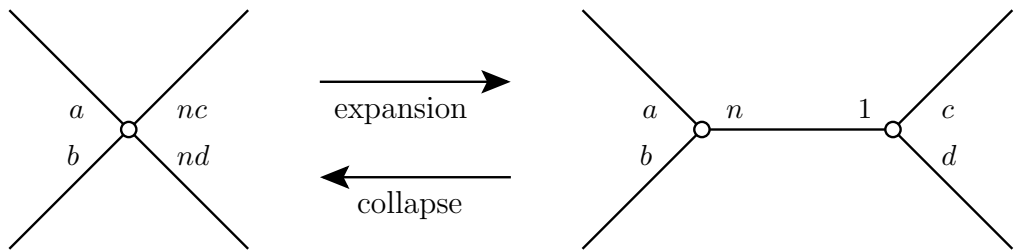}
\caption{The effect of elementary moves on
labeled graphs.}\label{fg:def-effect}
\end{figure}

\begin{figure}[ht]
\centering
\psfrag{slide}{slide}
\psfrag{induction}{induction}
\psfrag{A}{$\As$}
\psfrag{B}{$\As^{-1}$}
\psfrag{a}{{\small$a$}}
\psfrag{b}{{\small$b$}}
\psfrag{na}{{\small$na$}}
\psfrag{nb}{{\small$nb$}}
\psfrag{1}{{\small 1}}
\psfrag{lm}{{\small$\ell m$}}
\psfrag{ma}{{\small$ma$}}
\psfrag{k}{{\small$k$}}
\psfrag{m}{{\small$m$}}
\psfrag{klm}{{\small$k\ell m$}}
\includegraphics{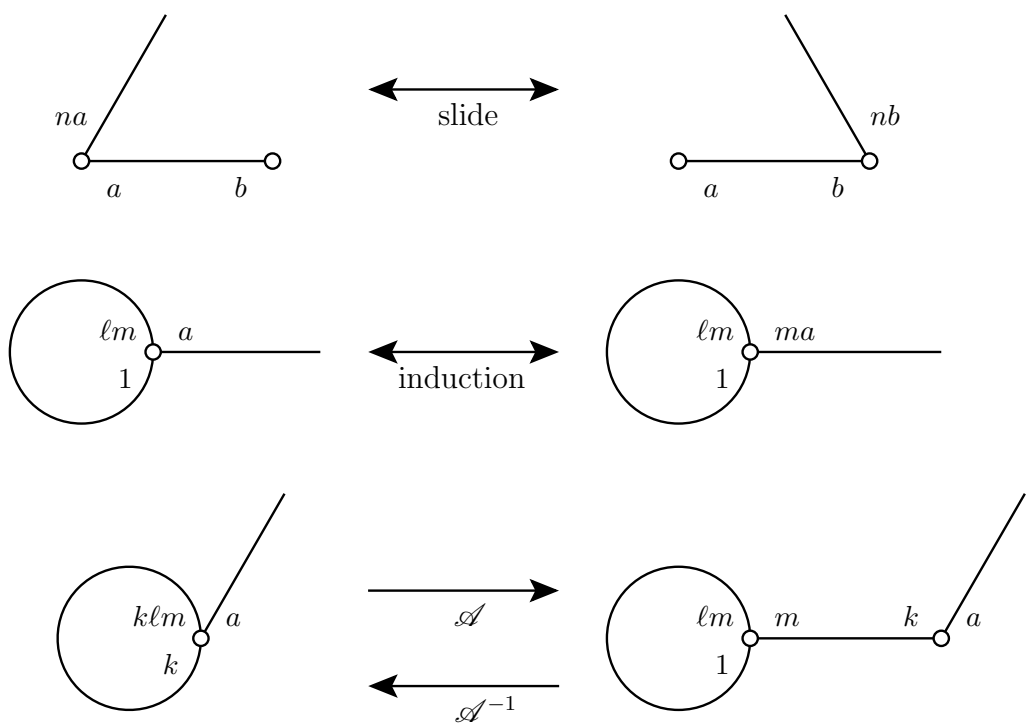}
\caption{The effect of the special elementary deformations described
  in Section \ref{ssc:pre} on labeled graphs.  These three moves
  suffice to relate any two reduced labeled graphs representing the
  same GBS group.}\label{fg:moves-effect}
\end{figure}

There is a certain bit of ambiguity in the function $\lambda$
resulting from different choices of generators for $G_{e}$ and
$G_{v}$.  The different choices result in changing the signs of
$\lambda(e)$ for all $e \in E_{o}(v)$ at some vertex $v$ or changing
the signs of $\lambda(e)$ and $\lambda(\bar{e})$ for some edge $e$.
Such changes are called \textit{admissible sign changes}.  We consider
two labeled graphs the same if they differ by admissible signs
changes.  Our labeled graphs are always equipped with a marking, i.e.,
there is $G$--tree $T$ with $T/G = (\Gamma,\lambda)$ as graphs of
groups.  We record in Figure \ref{fg:def-effect} and Figure
\ref{fg:moves-effect} the effect of the elementary moves and the
special moves listed in Section \ref{ssc:pre} on labeled graphs
\cite{ar:CFpp,ar:F06}.

A GBS group $G$ is \textit{non-elementary} if it is not isomorphic to
$\Z, \Z^{2}$ or the Klein bottle group.  For a non-elementary GBS
group, the elliptic subgroups arising from a labeled graph are
determined algebraically and do not depend on the particular tree
\cite{ar:F03}.  This implies that any two (marked) labeled graphs for
a GBS group are related by a sequence of elementary moves.  Further,
the set of conjugacy classes of elliptic subgroups is fixed by any
outer automorphism.  Thus if $\DS$ is the deformation space containing
one of these trees, then $\Out(G)_{\DS} = \Out(G)$.  When we speak of
a deformation space for a GBS group, we will always mean this
particular canonical deformation space.  For another class of groups
that have $\Out(G)_{\DS} = \Out(G)$ for a particular deformation space
$\DS$ see \cite{ar:C07}.

In the following example we build a star in $K$ and in $W$ and
describe the deformation retract $K \to W$ on these stars for the
group $G = BS(2,4)$.  Given a set $S$ of vertices in a simplicial
complex, the subcomplex \textit{spanned by S} is the subcomplex
consisting of all simplices whose vertices belong to the set $S$.

\begin{example}\label{ex:star}
  Let $G = BS(2,4)$ and $\DS$ be the canonical deformation space for
  $BS(2,4)$.  Fix $T \in W$ that has associated labeled graph as
  pictured on the left in Figure \ref{fg:label}.  Pick a vertex $v \in
  V(T)$ and choose a generator $g$ of $G_{v}$.  Label the edges
  emanating from $v$ as $E_{o}(v) =
  \{e_{0},e_{1},f_{0},f_{1},f_{2},f_{3}\}$ where $ge_{s} =e_{s+1\mod
    2}$ and $gf_{s} = f_{s+1 \mod 4}$.

  Two expansions of $T$ are defined by using pairs $J_0 = (\la
  g^2\ra,\{e_0,f_0,f_2\})$ and $J_1 = (\la g^2\ra,\{e_0,f_1,f_3\})$.
  Expansion by either of these pairs results in the labeled graph in
  the bottom left of Figure \ref{fg:bs24}.  However, the marking is
  different.  See Figure \ref{fg:j0j1} for a local picture of the
  expansion at $v$ to see the difference.  A third expansion of $T$ is
  defined by the pair $I = (\la g^2\ra,\{f_0,f_2\})$.  A final
  expansion is given by the pair $I' = (\la g \ra,
  \{f_0,f_1,f_3,f_4\})$.  This is the same as the expansion using $(\la
  g \ra,\{e_0,e_1\})$. Each of the trees $T^{J_0}$, $T^{J_1}$ and
  $T^{I'}$ can be further expanded using $I$.  These resulting labeled
  graphs are pictured in Figure \ref{fg:bs24}.

\begin{figure}[ht]      
  \centering 
  \psfrag{A}{$T^{J_0,I},T^{J_1,I}$}
  \psfrag{B}{$T^I$}
  \psfrag{C}{$T^{I,I'}$}
  \psfrag{E}{$T^{J_0},T^{J_1}$}
  \psfrag{F}{$T$}
  \psfrag{G}{$T^{I'}$}
  \psfrag{1}{{\small 1}}
  \psfrag{2}{{\small 2}}
  \psfrag{4}{{\small 4}}
  \includegraphics{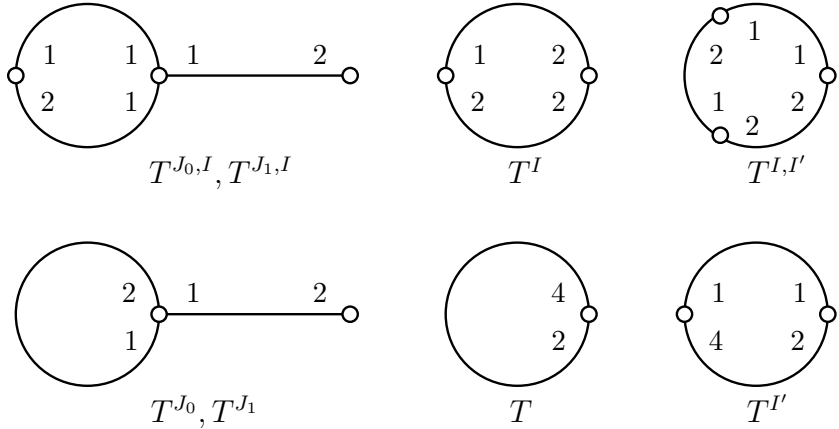}
  \caption{The labeled graphs for $BS(2,4)$ representing the vertices
    in the star of $T$ in $K$.}\label{fg:bs24}
\end{figure}

\begin{figure}[ht]
  \centering 
  \psfrag{a}{{\small$e_0$}} 
  \psfrag{b}{{\small$e_1$}}
  \psfrag{A}{$J_0$}
  \psfrag{B}{$J_1$} 
  \psfrag{c}{{\small$f_0$}} 
  \psfrag{d}{{\small$f_1$}}
  \psfrag{e}{{\small$f_2$}}
  \psfrag{f}{{\small$f_3$}}
  \includegraphics{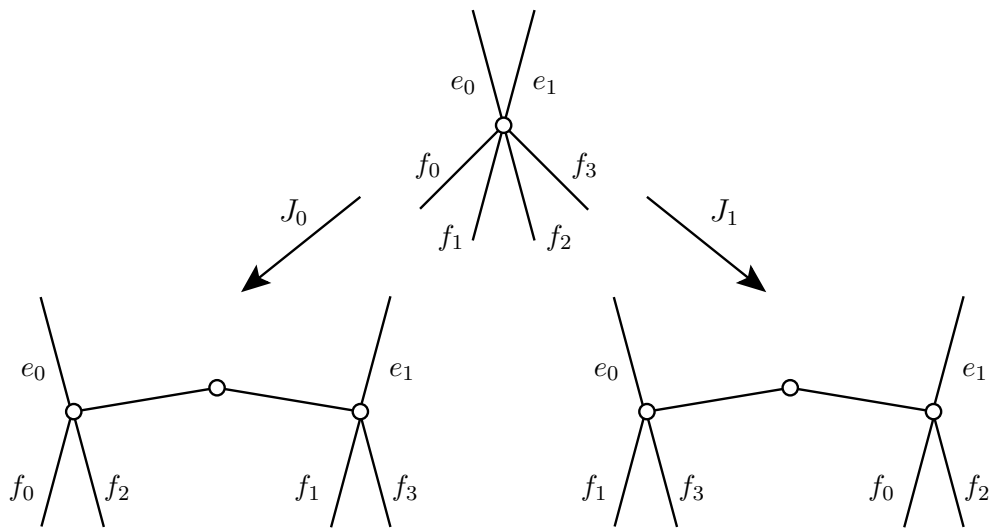}
  \caption{The expansions defined by $J_0$ and $J_1$.}\label{fg:j0j1}
\end{figure}

By examining these labeled graphs is it apparent that there are no
other expansions that do not create a subdivision vertex.  Therefore
the star of $T$ in $K$ is the complex pictured in Figure
\ref{fg:oneedgestar}.

% \begin{figure}[h]
% \centering
% \input{star.eepic} \caption{The star diagram showing compatibility for
% the pre-ideal edges $j_{0},j_{1}$ and $i$.  The circles enclose the
% $\sigma_{k}$--preimages of the cosets $\la x \ra / \la x^{2} \ra$ for
% $k = i,j_{0}$ and $j_{1}$.}\label{fg:star}
% \end{figure}

\begin{figure}[ht]
  \centering
  \psfrag{T}{$T$}
  \psfrag{A}{$T^{J_1}$}
  \psfrag{B}{$T^{J_0}$}
  \psfrag{C}{$T^{J_0,I}$}
  \psfrag{D}{$T^{J_1,I}$}
  \psfrag{E}{$T^{I,I'}$}
  \psfrag{F}{$T^{I'}$}
  \psfrag{G}{$T^I$}
  \includegraphics{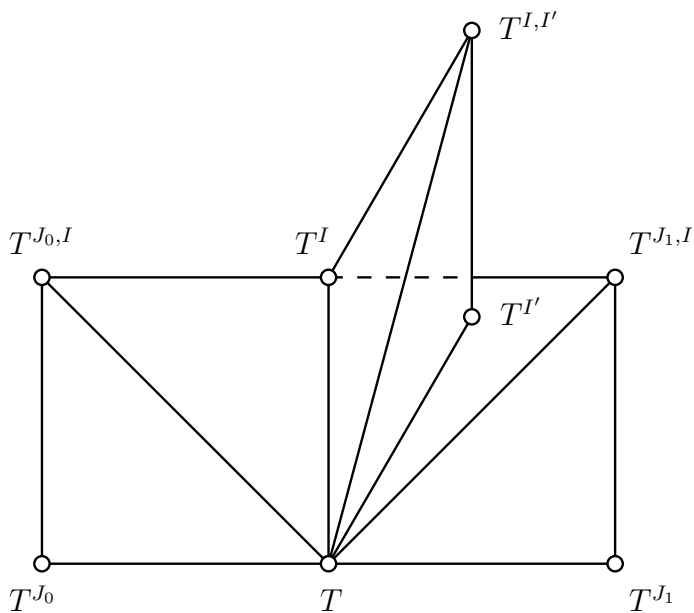} 
  \caption{The star of $T$ in $K$.  The star of $T$ in $W$ is the
    ``V'' subcomplex spanned by $T^{J_0,I}, T,
    T^{J_1,I}$.}\label{fg:oneedgestar}
\end{figure}

Of these labeled graphs, only those for $T$, $T^{J_0,I}$ and
$T^{J_1,I}$ are covered by shelters.  Therefore by Corollary
\ref{co:tw=t}, the star of $T$ in the complex $W$ is the two edge
``V'' subcomplex pictured Figure \ref{fg:oneedgestar} spanned by the
vertices labeled by $T^{J_0,I},T,T^{J_1,I}$.  The deformation retract
$h$ sends the vertices $T^{J_0},T^{I},T^{J_1},T^{I'}$ and $T^{I,I'}$
to the vertex $T$, fixing the other vertices.
\end{example}

By Theorem \ref{th:fd}(2) the deformation retract $W$ in Example
\ref{ex:star} is finite dimensional.  We remark that both $\DS$ and
the spine $K$ in this case are not finite dimensional.  We now present
an example of a deformation space $\DS$ for which the deformation
retract $W$ is not finite dimensional.

\begin{example}\label{ex:infinite-d}Let $G$ be the GBS group
  defined by $\la x,t_{0},t_{1} \ | \ t_{0}xt_{0}^{-1} =
  t_{1}xt_{1}^{-1} = x^{2} \ra$.  This presentation is represented by
  the labeled graph on the left in Figure \ref{fg:infinite-d}.
  Clearly the canonical deformation space $\DS$ containing this
  $G$--tree is neither non-ascending nor has $b_{1}(\DS) \leq 1$.  Let
  $T_{1}$ be a $G$--tree in $W$ with associated labeled graph pictured
  in the left in Figure \ref{fg:infinite-d}.  Pick a vertex $v \in
  V(T)$ and choose a generator $g$ of $G_v$.  Label the edges
  emanating from $v$ are labeled $E_{o}(v) =
  \{e,e_{0},e_{1},f,f_{0},f_{1}\}$ where $ge = e, gf = e, ge_{s} =
  e_{s+1 \mod 2}, gf_{s} = f_{s+1 \mod 2}$.  Inductively, let $T_k$ be
  the result of of sliding $e_0$ across $f$ in $T_{k-1}$.  The labeled
  graph for $T_k/G$ is the labeled graph in the center in Figure
  \ref{fg:infinite-d}.

  There is a collapse sequence $T^{k}_k \to \cdots \to T_k^\ell \to
  \cdots \to T_k^0 = T_k$ where the $G$--trees $T^\ell_k$ have
  associated labeled graphs as in the right in Figure
  \ref{fg:infinite-d} for $0 \leq \ell \leq k$.  Each of these
  $G$--trees is covered by shelters of type S3, and hence they
  represent vertices in $W$ that span a $k$--simplex.  Therefore $W$
  contains simplices of arbitrarily high dimension and is therefore not
  finite dimensional.

\begin{figure}[ht]
  \centering
  \psfrag{A}{$T_0$}
  \psfrag{B}{$T_k$}
  \psfrag{C}{$T_k^{\ell}$}
  \psfrag{1}{{\small 1}}
  \psfrag{2}{{\small 2}}
  \psfrag{2k}{{\small $2^k$}}
  \psfrag{2l}{{\small $2^{k-\ell}$}}
  \psfrag{o}{$\cdot$}
  \includegraphics{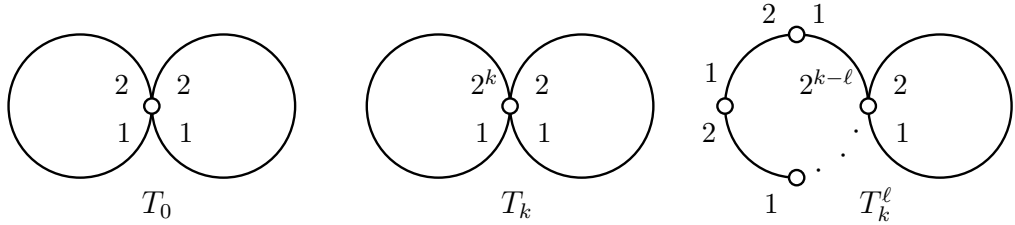}
  \caption{Labeled graphs representing the GBS group in Example
    \ref{ex:infinite-d}.  The deformation retract $W$ for the
    canonical deformation space for this group is not finite
    dimensional.}\label{fg:infinite-d}
\end{figure}
\end{example}

In the next section we will construct the deformation retract $W$ of
the canonical deformation space for $BS(p,q)$.  For another example
see \cite{ar:CT}, where the deformation retract $W$ of the canonical
deformation space for the GBS group $G = \la x,y,z \ | \ x^{n} = y^{n}
= z^{n} \ra$ for any $n \geq 2$ is constructed.

%%%%%%%%%%%%%%%%%%%%%%%%%%%%%%%%%%%%%%%%%%%%%%%%%%%%%%%%%%%%%%%%%%%%%%%%%%%%% 

\section{The canonical deformation space for
  $BS(p,q)$}\label{sc:dbspq}

For the remainder we use $G_{p,q}$ to denote the non-solvable
Baumslag--Solitar group $BS(p,q)$ where $|q| \geq p > 1$ (see
\eqref{eq:bs} in the Introduction).  Our aim is to compute the
automorphism groups $\Aut(G_{p,q})$ using the action on the canonical
deformation space associated to $G_{p,q}$.  We will denote this
deformation space by $\DS_{p,q}$.  The deformation retract for
$\DS_{p,q}$ described in the Section \ref{ssc:def} is denoted
$W_{p,q}$.  To begin we separate the non-solvable Baumslag--Solitar
groups into three types:

\begin{itemize}

\item[1.]  $p$ does not divide $q$; 

\item[2.] $p = |q|$; 

\item[3.] $q = pn$ where $|n| > 1$.

\end{itemize}

The first type was originally studied by Collins \cite{ar:C78} in the
special case that $p$ and $q$ do not share any factors.  This case was
analyzed fully by Gilbert, Howie, Metaftsis and Raptis
\cite{ar:GHMR00}.  Gilbert et al.~show that if $p$ does not divide $q$
then the group $\Aut(G_{p,q})$ acts on the same $G_{p,q}$--tree that
appears in Figure \ref{fg:label} \cite{ar:GHMR00}.  This observation
was independently discovered by Pettet \cite{ar:P99}.  Gilbert et
al.~analyze this action to show that $\Out(G_{p,q})$ is isomorphic to
the dihedral group of order $2|q - p|$.  A presentation for the
automorphism group follows from this.

In the language of deformation spaces, the theorem of Gilbert et
al.~and independently Pettet translates to showing that the
deformation space $\DS_{p,q}$ is rigid when $p$ does not divide $q$.
A deformation space $\DS$ is \textit{rigid} if there is a unique
reduced $G$--tree in $\DS$, up to equivariant homeomorphism.
Therefore, as $\Out(G)_{\DS}$ acts on the deformation space,
preserving the set of reduced $G$--trees, this unique reduced
$G$--tree is fixed by $\Out(G)_{\DS}$.  Therefore $\Out(G)_{\DS}$ acts
on the unique reduced $G$--tree, extending the action of $G$
\cite{ar:BJ96}.  In the case of a rigid deformation space, the complex
$W$ is a single point representing the unique reduced $G$--tree.  The
computation of Gilbert et al.~is translated as the computation of the
stabilizer of this point.  Although we do not need it in what follows,
we remark that Levitt has given a complete classification rigid
deformation spaces \cite{ar:L02} (see also \cite{ar:CFpp}).

Gilbert et al.~computed the automorphism group of the second type with
a similar computation as for the first type \cite{ar:GHMR00}.  Once
again, the deformation space $W_{p,q}$ is rigid \cite{ar:L02}.  Thus
the stabilizer of this unique reduced $G_{p,q}$--tree, is the entire
group $\Out(G_{p,q})$.  The computation of this stabilizer is similar
to the computation for the other stabilizers of other $G_{p,q}$--trees
that we show later on.  If $p = q$ then $\Out(G_{p,p}) = \Z \rtimes
(\Z_{2} \times \Z_{2})$, if $p = -q$ then $\Out(G_{p,-p}) = \Z_{2p}
\rtimes \Z_{2}$.  The appearance of the $\Z$ factor in $\Out(G_{p,p})$
is due to the fact that $G_{p,p}$ has a nontrivial center.  Computing
the full automorphism groups from here is trivial.

The third type was studied by Collins and Levin \cite{ar:CL83} and a
presentation for the group $\Aut(G_{p,q})$ was given by algebraic
means.  We will approach this using deformation spaces.  By passing to
an invariant tree $X_{p,q} \subseteq W_{p,q}$ ($X_{p,q} = W_{p,q}$ if
$|q|/p$ is prime) we will compute $\Out(G_{p,q})$ via Bass--Serre
theory.  The presentations of $\Out(G_{p,q})$ and $\Aut(G_{p,q})$ in
this case are given in Theorem \ref{th:presentations}.  For the
remainder we suppose that $q = pn$ where $p,|n| > 1$.

%%%%%%%%%%%%%%%%%%%%%%%%%%%%%%%

\subsection{$G_{p,q}$--trees in $W_{p,q}$}\label{ssc:Gpqtrees}

Using the results of Section \ref{ssc:lkW} we will give a
classification of the $G_{p,q}$--trees representing vertices in
$W_{p,q}$.  Denote by $\Pi(n) \subset \Z$ the multiplicative monoid
generated by the factors of $n$, and $\Pi^+(n) = \{ m \in \Pi(n) \, |
\, m \geq 1 \}$.

\begin{figure}[ht]
  \centering
  \psfrag{1}{{\small 1}}
  \psfrag{p}{{\small $p$}}
  \psfrag{q}{{\small $q$}}
  \psfrag{n}{{\small $n$}}
  \psfrag{m}{{\small $m$}}
  \includegraphics{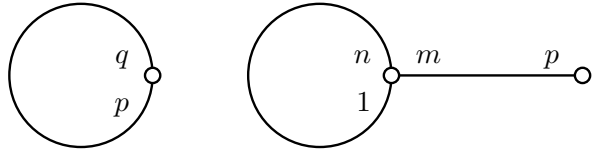}
  \caption{Labeled graphs representing reduced $G_{p,q}$--trees in
    $W_{p,q}$, $m \in \Pi^+(n), m \neq 1$.}
  \label{fig:reduced}
\end{figure}

\begin{lemma}\label{lm:reduced}
  If $T \in W_{p,q}$ is reduced, then after some admissible sign
  changes the associated labeled graph for $T/G_{p,q}$ is either:

  \begin{itemize}
  \item[1.] a single edge $e$ with $\iv(e) = \tv(e)$ and labels
    $\lambda(e) = p,\lambda(\bar{e}) = q$; or

  \item[2.] two edges $e,f$ with $\iv(e) = \tv(e) = \iv(f)$ and labels
    $\lambda(e) = 1,\lambda(\bar{e}) = n,\lambda(f) = m \neq 1$ where
    $m \in \Pi^+(n)$ and $\lambda(\bar{f}) = p$.
  \end{itemize}
  See Figure \ref{fig:reduced}.
\end{lemma}

\begin{proof}
  We will show that the collection of such labeled graphs is closed
  under slides, inductions and $\As^{\pm 1}$--moves.  By Theorem
  \ref{th:ais}, this implies the conclusion of the lemma.

  If $T$ is as in case 1, the only possible move is an $\As$--move.
  This results in a labeled graph as described in case 2, except
  possibly $\lambda(f)$ could be negative.  If this factor is
  negative, changing the signs of the labels on the edge $f$ and then
  changing the sign of $\lambda(\bar{f})$ results in a labeled graph
  as in case 2.  

  If $T$ is as in case 2 then a slide move possibly followed by
  admissible sign changes as above either results in a labeled graph
  as in case 2 with $\lambda(f) = |nm|$ or $|m/n|$.  The latter case
  is only possible when $n$ properly divides $\lambda(f)$ as otherwise
  the resulting labeled graph is not reduced.  An induction move only
  changes the label $\lambda(f)$, which is either multiplied or
  divided by a factor of $n$.  The latter case is only possible when
  this factor is not equal to $\lambda(f)$ as otherwise the resulting
  labeled graph is not reduced.  Thus again after possible admissible
  sign changes the resulting labeled graph is also as in case 2.  An
  $\As^{-1}$--move is only possible when $\lambda(f) \leq |n|$.  The
  resulting labeled graph is as in case 1.  An $\As$--move is not
  possible in this case.
\end{proof}

Given such labeled graphs as in case 2 of Lemma \ref{lm:reduced}, we
say an induction move is \textit{increasing} if the label
$|\lambda(f)|$ is larger after the induction move and
\textit{decreasing} otherwise.  Similarly define \textit{increasing
  slides} and \textit{decreasing slides}.

\begin{lemma}\label{lm:Gpqtrees}
  If $T \in W_{p,q}$, then after some admissible sign changes the
  associated labeled graph for $T/G_{p,q}$ is either as in case 1 of
  Lemma \ref{lm:reduced} or it consists of a shelter $\gamma$ of type
  S3 with a single edge $f$ attached at $\iv(f)$.  The labels on $f$
  are $\lambda(f) = m$ where $m \in \Pi^+(n)$ and
  $\lambda(\bar{f})=p$.
\end{lemma}

\begin{figure}[ht]
\centering
\psfrag{n1}{{\small $n_1$}}
\psfrag{n2}{{\small $n_2$}}
\psfrag{nm}{{\small $n_\ell$}}
\psfrag{p}{{\small $p$}}
\psfrag{1}{{\small $1$}}
\psfrag{m}{{\small $m$}}
\includegraphics{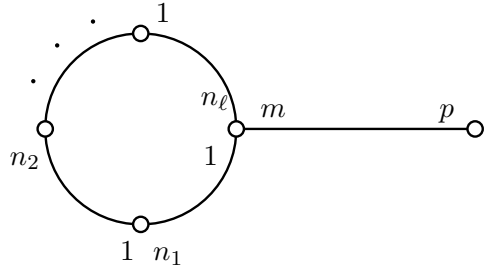} 
\caption{Labeled graphs representing $G_{p,q}$--trees in $W_{p,q}$.
  The labels on the circuit $\gamma$ satisfy $n_1\cdots n_\ell = n$
  and $|n_s| \neq 1$ with the possible exception of
  $n_\ell$.}\label{fg:linkT}
\end{figure}

\begin{proof}
  First, suppose the labeled graph representing $T$ is a cycle
  $\gamma$.  If $\gamma$ is covered by shelters of type S1, then
  $\gamma$ collapses to a cycle with at least two edges, contradicting
  Lemma \ref{lm:reduced}.  If $\gamma$ is a shelter of type S3 then
  after collapsing every edge except one, we get a reduced labeled
  graph that is an ascending loop, contradicting Lemma
  \ref{lm:reduced}.  If $\gamma$ is a shelter of type S2, then as
  there are no subdivision vertices in $G_{p,q}$--trees in $W_{p,q}$,
  hence $\gamma$ is a single edge and hence must be as in case 1 of
  Lemma \ref{lm:reduced}.

  Otherwise, the labeled graph consists of a circuit $\gamma$ with
  some finite trees attached.  Since there is a unique conjugacy class
  of maximal elliptic subgroups for $G_{p,q}$, there is at most one
  valence one vertex in a labeled graph representing a $G_{p,q}$--tree
  in $W_{p,q}$.  Therefore, there is at most one finite tree, $F$,
  attached to $\gamma$ and it is linear.  As there is a unique maximal
  conjugacy class of elliptic subgroups and as $G_{p,q}$--trees in
  $W_{p,q}$ do not subdivision vertices, $F$ cannot have any valence
  two vertices and is hence a single edge $f$.

  First we suppose $\gamma$ contains a shelter of type S1.  Then
  $\gamma$ is covered by disjoint shelters of type S1, further there
  are at least two shelters of type S1 needed to cover $\gamma$.
  Therefore, we can reduce the labeled graph for $T$ to a get a
  reduced labeled graph with a circuit with at least two edges.  By
  Lemma \ref{lm:reduced}, this is a contradiction.  Therefore,
  $\gamma$ is either a type S2 or type S3 shelter.

  Now suppose $\gamma$ is a type S2 shelter.  Since there are no
  subdivision vertices and there is a single edge attached, $\gamma$
  is either one or two edges.  If $\gamma$ is a single edge, then as
  $f$ must be in a shelter and using Lemma \ref{lm:reduced} we see
  that the labels on $\gamma$ are $1$ and $n$.  Hence $\gamma$ can
  also be thought of as a shelter of type S3.  If $\gamma$ is two
  edges then a similar argument shows $\gamma$ can be thought of as a
  type S3 shelter.  Therefore we have shown that the labeled graph is
  shelter of type S3 with a single edge attached.

  To see which labels can appear on $f$ we collapse the cycle $\gamma$
  to a single edge.  The labels appearing on $\gamma$ are all factors
  of $n$, therefore collapsing $\gamma$ cannot change whether or not
  $\lambda(f) \in \Pi(n)$.  Since by Lemma \ref{lm:reduced},
  $\lambda(f) \in \Pi^+(n)$ after collapsing $\gamma$, we must have
  that $f \in \Pi^+(n)$ initially.  Also collapsing $\gamma$ does not
  change $\lambda(\bar{f})$.  Again, by Lemma \ref{lm:reduced}
  $\lambda(\bar{f}) = p$ initially.
\end{proof}

%%%%%%%%%%%%%%%%%%%%%%%%%%%%%%%

\subsection{The complex $X_{p,q}$}\label{ssc:Xpq}

\begin{figure}[ht]
  \psfrag{n}{{\small $n$}}
  \psfrag{nk}{{\small $|n^k|$}}
  \psfrag{nk1}{{\small $|n^{k-1}|$}}
  \psfrag{1}{{\small 1}}
  \psfrag{p}{{\small $p$}}
  \psfrag{q}{{\small $q$}}
  \centering
  \includegraphics{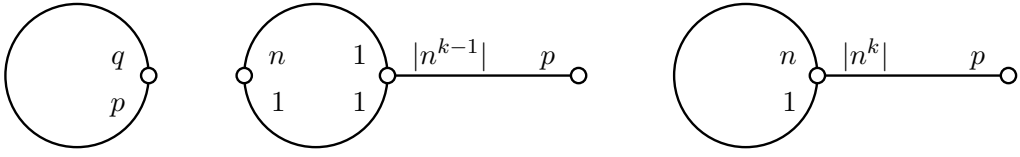}
  \caption{Labeled graphs representing $G_{p,q}$--trees in $W_{p,q}$
    when $|n|$ is prime ($k \geq 1$).}
  \label{fig:Xpqtrees}
\end{figure}

To see the motivation behind the definition of $X_{p,q}$ we will first
describe the complex $W_{p,q}$ when $|n|$ is prime.  As $|n|$ is
prime, by Lemma \ref{lm:Gpqtrees}, any $G_{p,q}$--tree in $W_{p,q}$
has associated labeled graph as pictured in Figure \ref{fig:Xpqtrees}.
We continue to use the notation from Lemma \ref{lm:reduced} to denote
the edges in the reduced labeled graphs.

As the longest collapse sequence between the labeled graphs in Figure
\ref{fig:Xpqtrees} is one, the complex $W_{p,q}$ is one dimensional,
hence a tree.  To each of the reduced labeled graphs in Figure
\ref{fig:Xpqtrees} we assign a non-negative integer called the
\textit{level}.  For the one edge reduced labeled graph, the level is
0.  For the two edge reduced labeled graphs, the level is the
non-negative integer $k$ appearing in the exponent of $n$.  We will
also use the term level when talking about a reduced $G_{p,q}$--tree
$T$ as the level of the associated labeled graph $T/G_{p,q}$.

As in Example \ref{ex:star}, the star of a reduced $G_{p,q}$--tree
with level 0 is a $p$--pod.  This is not the case when $|n|$ is
composite, although we will see in Lemma \ref{lm:disconnected} that in
this case the star is very close to being a $p$--pod.  The terminal
vertices of the $p$--pods are represented by labeled graphs pictured
in the center of Figure \ref{fig:Xpqtrees} with $k = 1$.  In addition
to collapsing to a $G_{p,q}$--tree with level 0, these trees also
collapse to a $G_{p,q}$--tree with level 1.

\begin{figure}[ht]
  \centering
  \psfrag{1}{{\small *}}
  \psfrag{2}{{\small $\phi_1$}}
  \psfrag{3}{{\small $\phi_2$}}
  \psfrag{4}{{\small $\phi_1\phi_2$}}
  \psfrag{5}{{\small $\phi_3$}}
  \psfrag{6}{{\small $\phi_1\phi_3$}}
  \psfrag{7}{{\small $\phi_2\phi_3$}}
  \psfrag{8}{{\small $\phi_1\phi_2\phi_3$}}
  \psfrag{l0}{{\small level 0}}
  \psfrag{l1}{{\small level 1}}
  \psfrag{l2}{{\small level 2}}
  \psfrag{l3}{{\small level 3}}
  \psfrag{T0}{{\small $T_0$}}
  \psfrag{T1}{{\small $T_1$}}
  \psfrag{T2}{{\small $T_2$}}
  \psfrag{T3}{{\small $T_3$}}
  \psfrag{c}{{\small $\cdot$}}
  \includegraphics{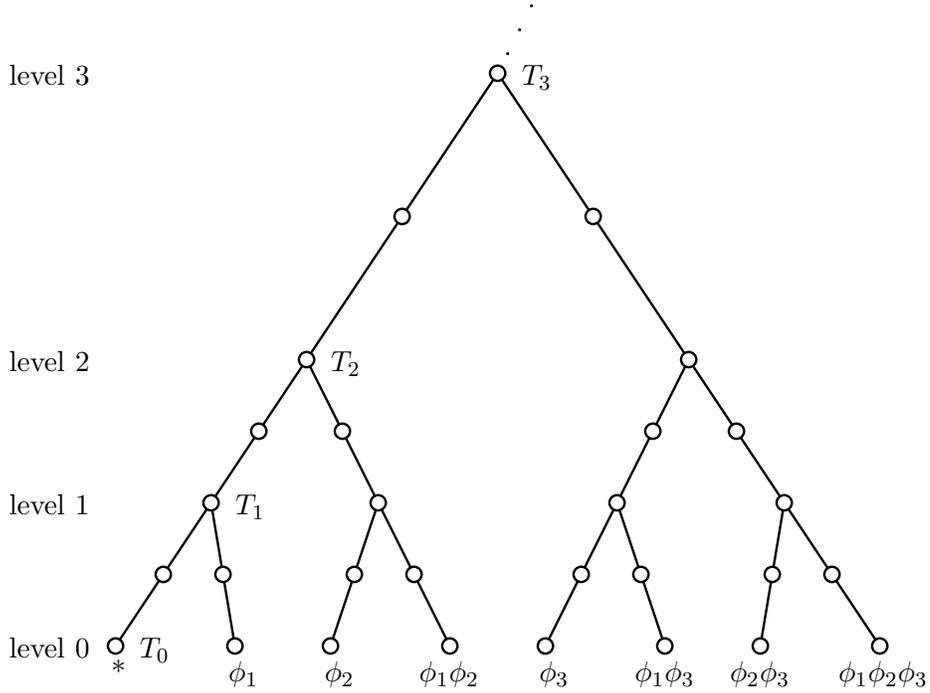}
  \caption{A piece of the complex $W_{2,4}$. The labels on the bottom
    show the action of the automorphisms $\phi_k$ listed in
    \eqref{eq:G_k}.}\label{fg:bs24sheet}
\end{figure}

The star of a reduced $G_{p,q}$--tree with level $k \geq 1$ is a
$(|n|+1)$--pod.  This follows as there are $|n|$ decreasing slides of
the edge $f$ counterclockwise around the loop $\bar{e}$ to a reduced
tree with level $k-1$ (further collapsing resulting in an
$\As^{-1}$--move is needed if $k=1$) and a unique increasing slide of
the edge $f$ clockwise around the loop $e$ to a reduced tree with
level $k+1$.  Figure \ref{fg:bs24sheet} shows a piece of the complex
$W_{2,4}$.  The entire complex $W_{2,4}$ is comprised of similar
pieces, glued $2$ at a time along a level 0 vertex, such that the
resulting complex is contractible.

We define $X_{p,q}$ to mimic $W_{p,q}$ in the case that $|n|$ is
prime.

\begin{definition}\label{def:X}
  Suppose $q = pn$ where $p,|n| > 1$.  Let $X_{p,q}$ be the subcomplex
  of $W_{p,q}$ spanned by $T \in W_{p,q}$ whose associated labeled
  graph is one of the three pictured in Figure \ref{fig:Xpqtrees}, $k
  \geq 1$.
\end{definition}

It is clear that $X_{p,q} = W_{p,q}$ when $|q|/p$ is prime.  In any
case, the longest collapse sequence between such labeled graphs is
one.  Hence the subcomplex $X_{p,q}$ is a (possibly disconnected)
graph.  The notion of level defined for reduced labeled graphs in
$W_{p,q}$ when $|n|$ is prime extends to reduced labeled graphs in
$X_{p,q}$.  We will show that $X_{p,q}$ is a tree.  First we show that
$X_{p,q}$ is connected.

\begin{lemma}\label{lm:klinduction}
  If $a,b$ divide $n$ and the product $ab$ also divides $n$ then the
  composition of the increasing inductions using the subgroups of
  index $a$ and $b$ is the increasing induction using the subgroup of
  index $ab$.  Further, if $ab = n$, then the composition of the
  increasing inductions is an increasing slide.
\end{lemma}

\begin{figure}[ht]
  \centering
  \psfrag{n}{{\small $n$}} 
  \psfrag{m}{{\small $m$}} 
  \psfrag{k}{{\small $a$}}
  \psfrag{l}{{\small $b$}} 
  \psfrag{kl}{{\small $ab$}}
  \psfrag{klm}{{\small $abm$}}
  \psfrag{km}{{\small $am$}}
  \psfrag{kn'}{{\small $a\ell$}}
  \psfrag{ln'}{{\small $b\ell$}}
  \psfrag{n'}{{\small $\ell$}}
  \includegraphics[width=15cm]{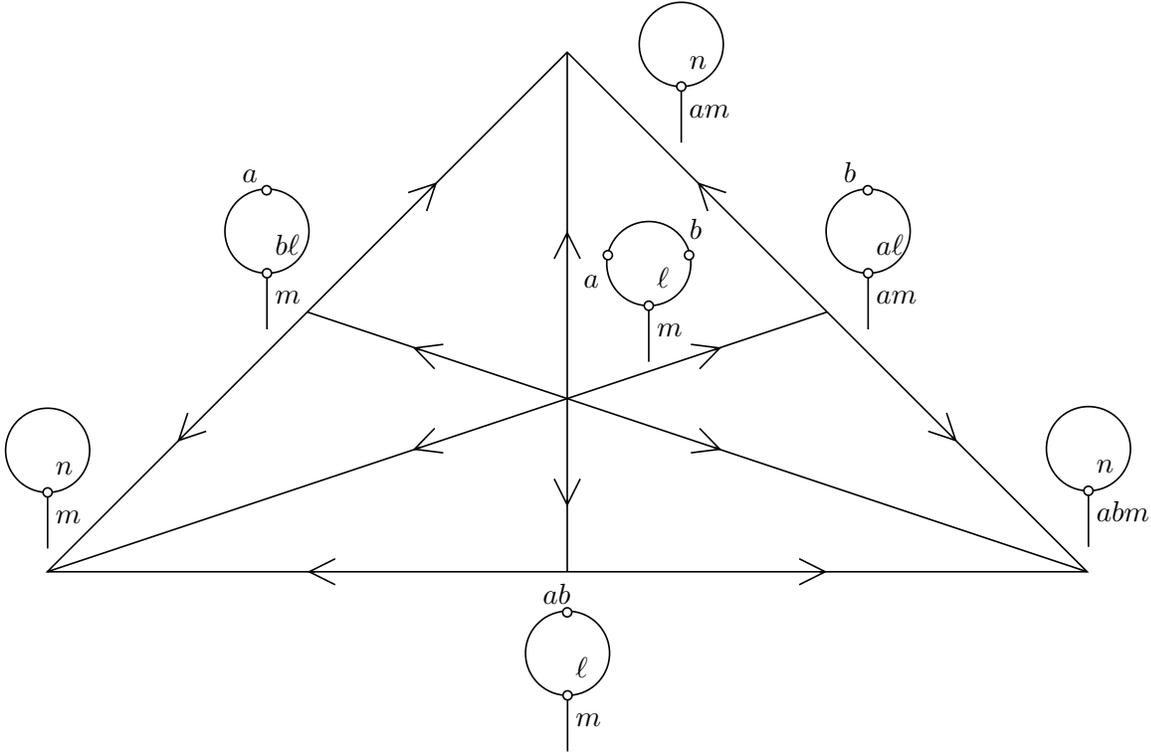}
  \caption{The 2--cell in $W_{p,q}$ in Lemma \ref{lm:klinduction}.
    The labels $1$ and $p$ are omitted from the
    figure.}\label{fig:2cell}
\end{figure}

\begin{proof}
  Write $n = ab\ell$.  The expansion appearing in an increasing
  induction depends only on the subgroup of the vertex group.  As
  there is a unique subgroup of any given index in $\Z$, the expansion
  is uniquely defined by the expansion in the labeled graph.  The
  2--cell in $W_{p,q}$ pictured in Figure \ref{fig:2cell} shows that
  the result of the composition of increasing inductions using the
  individual subgroups of index $a$ and $b$ is the same as the result
  of the increasing induction using the the subgroup of index $ab$.
  If $\ell = 1$, then the sequence of moves displayed along the bottom
  is increasing slide.
\end{proof}

\begin{lemma}\label{lm:increasinginductions}
  Let $T \in W_{p,q}$ be reduced and suppose $T' \in W_{p,q}$ is
  obtained by the composition of the increasing inductions using the
  subgroups of index $a$ and then $b$ and $T" \in W_{p,q}$ is obtained
  by the composition of the increasing inductions using the subgroups
  of index $b$ and then $a$.  Then $T' = T''$.
\end{lemma}

\begin{proof}
  If the product $ab$ divides $n$, then either composition of the
  increasing inductions is the increasing induction using the subgroup
  of index $ab$ by Lemma \ref{lm:klinduction}.  Hence the Lemma
  follows in this case.

  For the other case write $a = a'\ell$, $b = b' \ell$ where $\ell =
  \gcd(a,b)$.  Then $a'b'$ divides $n$.  The composition of the
  increasing inductions using $a$ and then $b$ can be written as the
  composition of the four increasing inductions using $\ell$, $a'$,
  $b'$ and then $\ell$.  Applying the earlier observation to the pair
  $a',b'$ we can write the composition of the increasing inductions as
  the the composition of the four increasing inductions using $\ell$,
  $b'$, $a'$ and then $\ell$.  But since these combine to give us the
  composition of the increasing inductions using $b'\ell = b$ and
  $a'\ell =a$ the Lemma follows.
\end{proof}

Define a map on reduced $G_{p,q}$--trees in $W_{p,q}$ to reduced
$G_{p,q}$--trees in $X_{p,q}$ by letting $x(T)$ be the tree obtained
by a composition of increasing inductions resulting in changing
$\lambda(f) = m$ to $mm' = n^k$ where $n$ does not divide $m'$.  By
Lemma \ref{lm:increasinginductions} $x$ is well-defined.  Denote the
integer $m'$ by $\lambda(T)$.

\begin{lemma}\label{lm:commute}
  If $T,T' \in W_{p,q}$ are related by a slide or an induction then
  $x(T)= x(T')$ or $x(T)$ and $x(T')$ are related by a slide.
\end{lemma}

\begin{proof}
  As by Lemma \ref{lm:klinduction} an increasing slide is a special
  case of an increasing induction we only need to show the lemma when
  $T'$ is obtained from $T$ by an increasing induction.

  If $\ell$ divides $\lambda(T)$, then $\lambda(T') = \lambda(T)/\ell$
  and by Lemma \ref{lm:increasinginductions} $x(T) = x(T')$ and the
  lemma holds.

  If $\ell$ does not divide $\lambda(T)$, then $\lambda(T') =
  n\lambda(T)/\ell$.  By Lemma \ref{lm:increasinginductions} the
  composition of the increasing inductions $T \to T' \to x(T')$ using
  the factors $\ell$ and then $\lambda(T')$ can be written as the
  composition of the increasing inductions $T \to x(T) \to x(T')$
  using the factors $\lambda(T)$ and then $n$.  Thus by Lemma
  \ref{lm:klinduction}, $x(T)$ and $x(T')$ are related by a slide.  
\end{proof}

\begin{proposition}\label{lm:connected}
The graph $X_{p,q}$ is connected.
\end{proposition}

\begin{proof}
  Let $T$ and $T'$ be reduced $G_{p,q}$--trees in $X_{p,q}$.  Then
  there is a 1--skeleton path in $W_{p,q}$ that connects these two
  trees.  Moreover, by Theorem \ref{th:ais} we can assume that this
  path is given by a sequence of slides, inductions and $\As^{\pm
    1}$--moves between reduced trees.

  By inserting an increasing induction move before any
  $\As^{-1}$--move and after any $\As$--move, we can assume that any
  $\As^{\pm 1}$--move appearing is between reduced trees in $X_{p,q}$.
  Therefore we only need to show that any two reduced $G$--trees in
  $X_{p,q}$ related by a sequence of slide and inductions between
  reduced $G_{p,q}$--trees in $W_{p,q}$ are in related by a sequence
  of slides between reduced $G_{p,q}$--trees in $X_{p,q}$.  By Lemma
  \ref{lm:commute} the map $x$ transforms such a sequence of
  inductions and slides into a sequence of slides between reduced
  $G_{p,q}$--trees in $X_{p,q}$.
\end{proof}

%The level of a reduced tree in $X_{p,q}$ is defined as before in the
%case that $|n|$ is prime.  
Therefore $X_{p,q}$ is a connected graph.  The $G_{p,q}$--trees in
$X_{p,q}$ form pieces in $X_{p,q}$ similar to the one pictured in
Figure \ref{fg:bs24sheet}.  This is true as from any reduced
$G_{p,q}$--tree with level $k \geq 1$ there are $|n|$ decreasing
slides resulting in a reduced $G_{p,q}$--tree with level $k-1$
(further collapsing resulting in an $\As^{-1}$--move is needed if
$k=1$) and a unique increasing slide resulting in a reduced
$G_{p,q}$--tree with level $k+1$.  In the next lemma we look at the
link of a $G_{p,q}$--tree with level 0.
    
\begin{lemma}\label{lm:disconnected}
  If $T \in W_{p,q}$ is a reduced with level $0$ then the link of $T$
  in $W_{p,q}$ is homotopy equivalent to a set of $p$ points.
  Moreover, each one of these points is naturally identified to an
  adjacent $G_{p,q}$--tree in $X_{p,q}$.
\end{lemma}

\begin{proof}
  By Lemma \ref{lm:Gpqtrees}, labeled graphs for $G_{p,q}$--trees in
  the link of $T$ all consist of a shelter of type S3 attached to a
  single edge $f$ at some vertex $\iv(f) = v$.  The labels on $f$ are
  $\lambda(f) = 1$ and $\lambda(\bar{f}) = p$.  Starting from $v$ the
  labels on the shelter are $(1,n_{1},1,n_{2},\ldots,1,n_{\ell})$
  where $n_{1}n_{2}\cdots n_{\ell} = n$ and $|n_{s}| \neq 1$ with the
  possible exception of $n_{\ell}$.  See Figure \ref{fg:linkT}.

  We will now define a deformation retraction of the link using
  Quillen's Poset Lemma (Lemma \ref{lm:poset}).  This retraction is
  the composition of two maps.  If $|n_\ell| \neq 1$, there is an
  elementary move that expands this $G_{p,q}$--tree $T'$ to a
  $G_{p,q}$--tree $h_{0}(T')$ where $n_{m} = 1$.  If $|n_\ell| = 1$,
  define $h_{0}(T') = T'$.  Then $h_{0}$ defines a poset map on the
  link that satisfies $h_{0}(T') \geq T'$.  Thus Quillen's Poset Lemma
  applies and defines a deformation retraction of the link.

  For $G_{p,q}$--trees $T'$ in the image of $h_{0}$ define $h_{1}(T')$
  as the $G_{p,q}$--tree found by collapsing the edges in the shelter
  of type S3 that are not adjacent to $v$.  Then $h_{1}$ defines a
  poset map on the image of $h_{0}$ satisfying $h_{1}(T') \leq T'$.
  Therefore Quillen's Poset Lemma applies again.  The image of
  $h_{1}h_{0}$ are the $p$ $G_{p,q}$--trees with labeled graph a shown
  in the middle of Figure \ref{fig:Xpqtrees} where $k=1$.  This proves
  the lemma.
\end{proof}

Using this lemma, we can prove that $X_{p,q}$ is a simply-connected,
hence a tree.

\begin{theorem}\label{th:Xtree}
  The subcomplex $X_{p,q} \subseteq W_{p,q}$ is a tree.
\end{theorem}

\begin{proof}
  As $X_{p,q}$ is a connected graph, we just need to show that
  $X_{p,q}$ does not contain a circuit, i.e.~an edge path
  homeomorphic to a circle.  Suppose there is a circuit $\gamma
  \subset X_{p,q}$.

  Since $W_{p,q}$ is simply-connected (Theorem \ref{th:cont}),
  $\gamma$ bounds a disk $D$ in $W_{p,q}$.  By Lemma
  \ref{lm:disconnected} and since $\gamma$ is homeomorphic to a
  circle, $\gamma$ cannot contain a vertex corresponding to a
  $G_{p,q}$--tree of level 0.  Therefore, $\gamma$ crosses some
  $G_{p,q}$--tree $T \in X_{p,q}$ with level $k \geq 1$ which is
  minimal among all $G_{p,q}$--tree along $\gamma$.  But as there is a
  unique edge from a $G_{p,q}$--tree with level $k$ to a
  $G_{p,q}$--tree with level $k+1$, this implies that the two edges of
  $\gamma$ adjacent to this tree are the same and therefore $\gamma$
  is not a circuit.  Hence $X_{p,q}$ does not contains a circuit.
\end{proof}

\begin{question}\label{qu:smallest}
  Is there a poset map that defines a deformation retraction from
  $W_{p,q} \to X_{p,q}$?  If so, can one define this map for a general
  deformation space to get a further deformation retraction of \DS?
%   If not, can the definition of $X_{p,q}$ be mimicked for a general
%   deformation space to get a contractible subcomplex that does not
%   involve induction moves?  Is such a subcomplex minimal?
\end{question}

%%%%%%%%%%%%%%%%%%%%%%%%%%%%%%%%%%%%%%%%%%%%%%%%%%%%%%%%%%%%%%%%%%%%%%%%%%%%% 

\section{Computation of ${\rm Out}(BS(p,q))$}\label{sc:out}

In this section we will use the action of $\Out(G_{p,q})$ on the tree
$X_{p,q}$ to give a presentation of this group in the case that $p$
properly divides $q$.  The vertices in $X_{p,q}$ corresponding to
$G_{p,q}$--trees that are not reduced are subdivision points, removing
these from $X_{p,q}$ does not alter the action of $\Out(G_{p,q})$.  We
begin by describing the quotient $X_{p,q}/ \Out(G_{p,q})$.

%%%%%%%%%%%%%%%%%%%%%%%%%%%%%%%

\subsection{The quotient $X_{p,q}/{\rm Out}(G_{p,q})$}\label{ssc:quot}

The following proposition is a restatement of a special case of
Proposition 5.3 in \cite{ar:T70} (cf.~\cite{ar:B93}).

\begin{proposition}\label{prop:quot}
  Suppose $G$ is a GBS group and $T,T'$ are $G$--trees with infinite
  cyclic point stabilizers.  If the associated labeled graphs for
  $T/G$ and $T'/G$ are isomorphic, then there is an outer automorphism
  $\Phi \in \Out(G)$ such that $T\Phi = T'$.
\end{proposition}
  
\begin{figure}[ht]
\centering
\psfrag{0}{{\small 0}}
\psfrag{1}{{\small 1}}
\psfrag{2}{{\small 2}}
\psfrag{3}{{\small 3}}
\psfrag{4}{{\small 4}}
\psfrag{5}{{\small 5}}
\psfrag{6}{{\small 6}}
\includegraphics{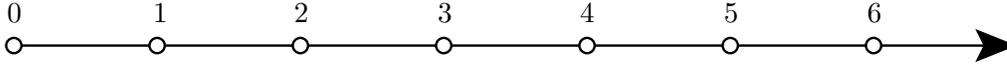}
\caption{The quotient graph
  $X_{p,q}/\Out(BS(p,q))$.}\label{fg:quotient}
\end{figure}

Therefore, the quotient $X_{p,q}/\Out(G_{p,q})$ can be identified with
the ray $[0,\infty)$, where the integer point $k$ is represented by an
(unmarked) labeled graph with level $k$.  See Figure
\ref{fg:quotient}.

%%%%%%%%%%%%%%%%%%%%%%%%%%%%%%%

\subsection{Infinite generation of ${\rm Out}(BS(p,q))$}\label{ssc:nfg}

Before we compute the vertex stabilizers of the tree $X_{p,q}$ we give
a geometric argument showing $\Out(G_{p,q})$ is not finitely generated
when $p$ properly divides $q$.  This follows easily from the following
lemma.

\begin{lemma}\label{lm:fg}
  Let $G$ be a finitely generated group acting by simplicial
  automorphisms on a connected simplicial complex $X$.  Then, for any
  point $x \in X$, there is a compact set $C \subset X$ containing $x$
  such that $GC$ is connected.
\end{lemma}

\begin{proof}
  This is Brown's finiteness criteria for type FP$_{0}$ \cite{ar:B87}.
  It is easy to prove in this case.  Without loss of generality, we
  can assume that $x \in X^{(0)}$.  Let $\{g_{0},\ldots,g_{m}\}$ be a
  finite generating set for $G$.  Then take $C$ to be the union of the
  1--skeleton paths between $x$ and $g_{s}x$ for $s = 0,\ldots,m$.
  Since the $g_{s}$ generate $G$, the set $GC$ is connected.
\end{proof}

Now to see that $\Out(G_{p,q})$ is not finitely generated, we apply
the above lemma to any $T \in X_{p,q}$ with level 0.  In any compact
set $C \subset X_{p,q}$ containing $T$ there is a $k \geq 0$ such that
any $G_{p,q}$--tree in $C$ has level at most $k$.  There is another
$T'$ with level 0 such that the geodesic from $T$ to $T'$ to passes
through a $G_{p,q}$--tree with level $k+1$.  By Proposition
\ref{prop:quot}, $T' \in \Out(G_{p,q})C$.  However, as $X_{p,q}$ is a
tree, any path from $T$ to $T'$ must go through a $G_{p,q}$--tree with
level $k+1$.  Such a path cannot lie entirely in $\Out(G_{p,q})C$ as
the action of $\Out(G_{p,q})$ preserves the level of a
$G_{p,q}$--tree.  In terms of Bestvina--Brady discrete Morse theory
\cite{ar:BB97}, the descending links of the $G_{p,q}$--trees with
level $k \geq 1$ are disconnected.  Hence we get the following
theorem, also noted by Collins and Levin \cite{ar:CL83}.

\begin{theorem}\label{th:nfg}
  If $p>1$ and $p$ properly divides $q$, then $\Out(BS(p,q))$ is not
  finitely generated.
\end{theorem}

%%%%%%%%%%%%%%%%%%%%%%%%%%%%%%%

\subsection{Vertex stabilizers in $X_{p,q}$}\label{ssc:stab}

The quotient $X_{p,q}/ \Out(G_{p,q})$ is a ray.  Lift this ray to a
ray in $X_{p,q}$.  Denote the $G_{p,q}$--tree on this ray representing
the integer point $k$ by $T_{k}$.  Without loss of generality, we can
assume $T_{0}/G_{p,q}$ gives rise to the presentation in
\eqref{eq:bs}.  Further, $T_k$ for $k \geq 1$ give rise to
presentations:
\begin{equation}\label{eq:bs2}
  G_{p,q} = \la a,b_k,t \, | \, a^p = b_k^{n^k}, tb_kt^{-1} = b_k^n \ra
\end{equation}
where $a \mapsto x, b_k\mapsto t^{-k}x^pt^k$ and $t\mapsto t$.

Let $H_k \subseteq \Out(G_{p,q})$ be the stabilizer of $T_{k}$.  We
have two cases, depending on whether $k=0$ or $k \geq 1$.  The
important fact we use is that if $\phi \in \Aut(G)$ fixes an
irreducible $G$--tree $T$ (here $G$ can be any finitely generated
group), then there is a unique $\phi$--equivariant simplicial
automorphism $f_{\phi}\co T \to T$ \cite{ar:B93}.  Thus we get an
action of the stabilizer $H_{T} \subseteq \Aut(G)$ on the tree $T$
that extends the action of $G$ (viewing $G/Z(G) = \Inn(G)$ as a
subgroup of $H_{T}$).

For an irreducible $G$--tree $T$, there are special automorphisms that
fix $T$.  These are called \textit{twists} as they generalize the
familiar notion of Dehn twist when the $G$--tree arises from a simple
closed curve on a surface.  We will only look a one type of twist, the
one that corresponds to a nonseparating curve, for a general
discussion of twists in $G$--trees, see \cite{ar:L05}.  Let $e$ be a
one edge loop in $T/G$ with vertex $v$ and stable letter $t$, when
viewing $e$ as giving rise to an HNN-extension.  Then for $z \in
G_{v}$ such that $zg = gz$ for all $g \in G_{e}$ the map that sends $t
\mapsto zt$ is a twist in $G$--tree $T$.  To fix some notation, for $g
\in G$ we denote the inner automorphism $g' \mapsto gg'g^{-1}$ by
$c_{g}$.
 
\medskip \noindent {\bf case 1}: $k = 0$.

We claim that $H_{0}$ is isomorphic to the dihedral group $\Z_{p|n-1|}
\rtimes \Z_{2}$, generated by the following automorphisms:
\begin{equation}\label{eq:G_0}
\begin{array}{rrlrrl}
\psi\co & x &\mapsto x  & \iota\co & x &\mapsto x^{-1} \\
        & t &\mapsto xt &          & t &\mapsto t
\end{array}
\end{equation}
Notice that $\psi^{p(n-1)} = c_{x}^{-p}$ and $\iota\psi =
\psi^{-1}\iota$.  Using normal forms for HNN-extensions, it is easy to
see that the outer automorphism class of $\psi^\ell$ for $1 \leq \ell
< p|n-1|$ is non-trivial.  Hence the image of $\la \psi,\iota \ra$ in
$\Out(G_{p,q})$ is the dihedral group $\Z_{p|n-1|} \rtimes \Z_{2}$.
The automorphism $\psi$ is a twist as described above for $T_{0}$,
hence it fixes $T_{0}$.  It is clear that the automorphism $\iota$
fixes any $G_{p,q}$--tree in $\DS_{p,q}$.  Thus the image of the
subgroup $\la \psi, \iota \ra$ is contained in $H_{0}$.  The claim
that this is an equality follows exactly as the computation of Gilbert
et al.~for the case when $p$ does not properly divide $q$
\cite{ar:GHMR00}.  This computation is similar to the computation in
the next case, thus we omit it.  \qed

\medskip \noindent {\bf case 2}: $k \geq 1$.

In this case, we claim that $H_k$ is isomorphic to the dihedral
group $\Z_{|n^k(n-1)|} \rtimes \Z_2$, generated by the following
automorphisms:
\begin{equation}\label{eq:G_k}
\begin{array}{rrlrrl} 
\phi_{k}\co & x &\mapsto x  & \iota\co & x &\mapsto x^{-1} \\
	& t &\mapsto (t^{-k}x^pt^k)t & & t &\mapsto t
\end{array}
\end{equation}
Notice that $\phi_{k+1}^n = \phi_k$ for $k \geq 1$, $\phi_1^n =
\psi^p$ and $\iota\phi_k = \phi_k^{-1}\iota$.  To prove this claim, it
is easier to use the presentations for $G_{p,q}$ in \eqref{eq:bs2}.
With this generating set, the automorphisms in \eqref{eq:G_k} are:
\begin{equation}\label{eq:G_k2}
\begin{array}{rrlrrl}
 & a &\mapsto a & & a &\mapsto a^{-1} \\
  \phi_{k}\co & b_{k} &\mapsto b_{k} & \iota\co & b & \mapsto
  b_{k}^{-1}\\
	& t &\mapsto b_{k}t & & t &\mapsto t
\end{array}
\end{equation}

Viewing these presentations as HNN-extensions $\la a,b_k \ra *_{\la
  tb_kt^{-1} = b_k^{n^k}\ra }$ and using normal forms for
HNN-extensions is easy to see that the outer automorphism class of
$\phi_k^\ell$ is non-trivial for $1 \leq \ell < |n^k(n-1)|$.  Hence
the image of $\la \phi_k,\iota \ra$ in $\Out(G_{p,q})$ is the dihedral
group $\Z_{|n^k(n-1)|} \rtimes \Z_{2}$.  Also with these
presentations, it is apparent that $\phi_{k}$ is a twist of the
$G_{p,q}$--tree $T_{k}$.  Thus the image of the subgroup $\la \phi_k,
\iota \ra$ is contained in $H_{k}$.  The action of the automorphisms
$\phi_{k}$ for $G_{2,4}$ on $X_{2,4}$ is shown in Figure
\ref{fg:bs24sheet}.

Now suppose that $\alpha \in H_{k}$.  Then lifting $\alpha$ to
$\Aut(G_{p,q})$, we have a $\alpha$--equivariant simplicial
automorphism $f_{\alpha}\co T_{k} \to T_{k}$.  There are two type of
vertices in $T_{k}$: those that are lifts of $\iv(f) = v$ or those
that are lifts of $\tv(f) = w$.  Lifts of $v$ belong to an axis of a
hyperbolic element of length one and lifts of $w$ do not, thus
$f_{\alpha}$ sends lifts of $v$ to lifts of $v$ and similar for $w$.
Therefore, after composing $\alpha$ with an inner automorphism, we can
assume that $f_{\alpha}$ fixes some lift of $v$ in $T$ with stabilizer
$\la b_k \ra$, which we continue to denote $v$.  Further, we can
assume that $v$ is adjacent to the unique vertex stabilized by $\la a
\ra$.  This implies that the axis of $t$ contains $v$.

We label the edges emanating from $v$ by $E_{o}(v) =
\{e,e_{0},\ldots,e_{|n|-1},f_{0},\ldots, f_{|n|^{k}-1} \}$ where $b_ke
= e,b_ke_s = e_{s+1 \mod |n|}$ and $b_kf_s = f_{s+1 \mod |n^k|}$.
Assume that $\tv(f_0)$ is stabilized by $\la a \ra$ and $tv =
\tv(e_0)$.  Define $\beta = c_{b_k} \phi_{k}^{n-1}$.  Then $\beta(a) =
b_{k}ab_{k}^{-1}$ and fixes $b_{k}$ and $t$.  After composing $\alpha$
with $\beta^{m'}$ for some $m'$ we can assume that $f_{\alpha}$ fixes
the edge $f_{0}$.  Hence, after composing with $\iota$ we can assume
that $\alpha(a) = a$ and $\alpha(b_{k}) = b_{k}$.

Therefore $f_{\alpha}$ permutes the edges $e_{s}$ for $s = 0,
\ldots,|n|-1$.  Now $\alpha(t)v = \alpha(t)f_{\alpha}(v) =
f_{\alpha}(tv) = \tv(e_{s}) = b_{k}^{s}tv$ for some $s$.  Thus
$t^{-1}b_{k}^{-s}\alpha(t) \in G_{v} = \la b_{k} \ra$.  Therefore
$t^{-1}b_{k}^{-s}\alpha(t) = b_{k}^{m}$.  Rewriting, we have
$\alpha(t) = b_{k}^stb_{k}^m = b_{k}^{s+mn}t = \phi_{k}^{s+mn}(t)$.
Thus $\alpha = \phi_{k}^{s+mn}$ and $H_{k}$ is as claimed.  \qed

Since $T_{k+1}$ is the unique $G_{p,q}$--tree of level $k+1$ adjacent
to $T_{k}$ we have inclusions $H_{k} \subseteq H_{k+1}$ for $k >0$.
Therefore, as a graph of groups, the infinite ray $X_{p,q} /
\Out(G_{p,q})$ collapses to a segment with one vertex corresponding to
$T_{0}$ and the other vertex corresponding to the end represented by
$(T_1,T_2,\ldots)$.  The stabilizer of this end is the direct limit:
\begin{equation*}\label{eq:dl}
    \lim_{\to} \Z_{|n^k(n-1)|} \rtimes \Z_{2} =
    \Z[\frac{1}{|n|}]/|n(n-1)| \rtimes \Z_{2}.
\end{equation*}
In the above, $1 \in \Z[\frac{1}{|n|}]/|n(n-1)|$ corresponds to the
outer automorphism class of $\phi_1$.
%%%%%%%%%%%%%%%%%%%%%%%%%%%%%%%

\subsection{Presentations}\label{ssc:presentations}

The computations from Sections \ref{ssc:quot} and \ref{ssc:stab} give
the presentation for $\Out(BS(p,q))$ appearing in the following
theorem.  The presentation for $\Aut(BS(p,q))$ follows routinely from
this.  This presentation was also found with an algebraic computation
by Collins and Levin \cite{ar:CL83}.

\begin{theorem}\label{th:presentations}
Let $q=pn$ where $p,|n| > 1$.  The automorphism group $\Aut(BS(p,q))$
is generated by the automorphisms $c_{x},c_{t},\psi,\iota,$ and
$\phi_{k}$ for $k \geq 1$ subject to the following relations:
\begin{equation*}\label{eq:aut}
\begin{array}{rl}
    c_{t}c_{x}^{p}c_{t}^{-1} = c_{x}^{q} & \iota^{-1} = \iota \\
    \iota c_{x} \iota = c_{x}^{-1} & \iota c_{t} \iota = c_{t} \\
    \iota \psi \iota = \psi^{-1} & \iota \phi_{k} \iota = 
    \phi_{k}^{-1} \mbox{, {\it for }} k \geq 1 \\
    \psi^{p} = \phi_{1}^{n} & \psi^{p(n-1)} = c_{x}^{-p} \\
    \psi c_{x} \psi^{-1} = c_{x} & \psi c_{t} \psi^{-1} = c_{x} c_{t} \\
    \phi_{k} c_{x} \phi_{k}^{-1} = c_{x} \mbox{, {\it for }} k \geq 1
    & \phi_{k} c_{t} \phi_{k}^{-1} = c_{t}^{-k} c_{x}^{p} c_{t}^{k+1}
    \mbox{, {\it for }} k \geq 1 \\
    \phi_{k+1}^{n} = \phi_{k} \mbox{, {\it for }} k \geq 1 &
\end{array}
\end{equation*}
The outer automorphism group has presentation:
\begin{equation*}\label{eq:out}
    \Out(BS(p,q)) = (\Z_{|p(n-1)|}*_{\Z_{|n-1|}} \Z[\frac{1}{|n|}]
    /|n(n-1)| \Z ) \rtimes \Z_{2}
  \end{equation*}
generated by the images of $\psi,\iota$ and $\phi_{k}$ for $k \geq 1$.
\end{theorem}

\begin{remark}\label{rm:opps}
  We remark that the relation $\iota \psi \iota = \psi^{-1}$ was
  omitted in Theorem 3.1 in \cite{ar:CL83}.  
\end{remark}

%%%%%%%%%%%%%%%%%%%%%%%%%%%%%%%%%%%%%%%%%%%%%%%%%%%%%%%%%%%%%%%%%%%%%%%%%%%%% 

\bibliographystyle{siam}
\bibliography{bibliography}

\def\cprime{$'$}
\begin{thebibliography}{10}

\bibitem{ar:B93}
{\sc H.~Bass}, {\em Covering theory for graphs of groups}, J. Pure Appl.
  Algebra, 89 (1993), pp.~3--47.

\bibitem{ar:BJ96}
{\sc H.~Bass and R.~Jiang}, {\em Automorphism groups of tree actions and of
  graphs of groups}, J. Pure Appl. Algebra, 112 (1996), pp.~109--155.

\bibitem{ar:BK90}
{\sc H.~Bass and R.~Kulkarni}, {\em Uniform tree lattices}, J. Amer. Math.
  Soc., 3 (1990), pp.~843--902.

\bibitem{ar:BS62}
{\sc G.~Baumslag and D.~Solitar}, {\em Some two-generator one-relator
  non-{H}opfian groups}, Bull. Amer. Math. Soc., 68 (1962), pp.~199--201.

\bibitem{ar:BB97}
{\sc M.~Bestvina and N.~Brady}, {\em Morse theory and finiteness properties of
  groups}, Invent. Math., 129 (1997), pp.~445--470.

\bibitem{ar:BF91}
{\sc M.~Bestvina and M.~Feighn}, {\em Bounding the complexity of simplicial
  group actions on trees}, Invent. Math., 103 (1991), pp.~449--469.

\bibitem{ar:B87}
{\sc K.~S. Brown}, {\em Finiteness properties of groups}, in Proceedings of the
  Northwestern conference on cohomology of groups (Evanston, Ill., 1985),
  vol.~44, 1987, pp.~45--75.

\bibitem{ar:C05}
{\sc M.~Clay}, {\em Contractibility of deformation spaces of {$G$}--trees},
  Algebr. Geom. Topol., 5 (2005), pp.~1481--1503 (electronic).

\bibitem{ar:CT}
\leavevmode\vrule height 2pt depth -1.6pt width 23pt, {\em Deformation spaces
  of {$G$}--trees}, PhD thesis, University of Utah, 2006.

\bibitem{ar:C07}
\leavevmode\vrule height 2pt depth -1.6pt width 23pt, {\em A fixed point
  theorem for deformation spaces of {$G$}-trees}, Comment. Math. Helv., 82
  (2007), pp.~237--246.

\bibitem{ar:CFpp}
{\sc M.~Clay and M.~Forester}, {\em Whitehead moves for {$G$}--trees}.
\newblock {\footnotesize\textbf{arXiv:0710.2107}}.

\bibitem{ar:C78}
{\sc D.~J. Collins}, {\em The automorphism towers of some one-relator groups},
  Proc. London Math. Soc. (3), 36 (1978), pp.~480--493.

\bibitem{ar:CL83}
{\sc D.~J. Collins and F.~Levin}, {\em Automorphisms and {H}opficity of certain
  {B}aumslag-{S}olitar groups}, Arch. Math. (Basel), 40 (1983), pp.~385--400.

\bibitem{ar:CM87}
{\sc M.~Culler and J.~W. Morgan}, {\em Group actions on {${\bf R}$}-trees},
  Proc. London Math. Soc. (3), 55 (1987), pp.~571--604.

\bibitem{ar:CV86}
{\sc M.~Culler and K.~Vogtmann}, {\em Moduli of graphs and automorphisms of
  free groups}, Invent. Math., 84 (1986), pp.~91--119.

\bibitem{ar:F02}
{\sc M.~Forester}, {\em Deformation and rigidity of simplicial group actions on
  trees}, Geom. Topol., 6 (2002), pp.~219--267 (electronic).

\bibitem{ar:F03}
\leavevmode\vrule height 2pt depth -1.6pt width 23pt, {\em On uniqueness of
  {JSJ} decompositions of finitely generated groups}, Comment. Math. Helv., 78
  (2003), pp.~740--751.

\bibitem{ar:F06}
\leavevmode\vrule height 2pt depth -1.6pt width 23pt, {\em Splittings of
  generalized {B}aumslag-{S}olitar groups}, Geom. Dedicata, 121 (2006),
  pp.~43--59.

\bibitem{ar:GHMR00}
{\sc N.~D. Gilbert, J.~Howie, V.~Metaftsis, and E.~Raptis}, {\em Tree actions
  of automorphism groups}, J. Group Theory, 3 (2000), pp.~213--223.

\bibitem{ar:GL072}
{\sc V.~Guirardel and G.~Levitt}, {\em Deformation spaces of trees}, Groups
  Geom. Dyn., 1 (2007), pp.~135--181.

\bibitem{ar:GL07}
\leavevmode\vrule height 2pt depth -1.6pt width 23pt, {\em The outer space of a
  free product}, Proc. Lond. Math. Soc. (3), 94 (2007), pp.~695--714.

\bibitem{ar:L05}
{\sc G.~Levitt}, {\em Automorphisms of hyperbolic groups and graphs of groups},
  Geom. Dedicata, 114 (2005), pp.~49--70.

\bibitem{ar:L02}
\leavevmode\vrule height 2pt depth -1.6pt width 23pt, {\em Characterizing rigid
  simplicial actions on trees}, in Geometric methods in group theory, vol.~372
  of Contemp. Math., Amer. Math. Soc., Providence, RI, 2005, pp.~27--33.

\bibitem{ar:L07}
\leavevmode\vrule height 2pt depth -1.6pt width 23pt, {\em On the automorphism
  group of generalized {B}aumslag-{S}olitar groups}, Geom. Topol., 11 (2007),
  pp.~473--515.

\bibitem{ar:MM96}
{\sc D.~McCullough and A.~Miller}, {\em Symmetric automorphisms of free
  products}, Mem. Amer. Math. Soc., 122 (1996), pp.~viii+97.

\bibitem{ar:P89}
{\sc F.~Paulin}, {\em The {G}romov topology on {${\bf R}$}-trees}, Topology
  Appl., 32 (1989), pp.~197--221.

\bibitem{ar:P99}
{\sc M.~R. Pettet}, {\em The automorphism group of a graph product of groups},
  Comm. Algebra, 27 (1999), pp.~4691--4708.

\bibitem{ar:Q78}
{\sc D.~Quillen}, {\em Homotopy properties of the poset of nontrivial
  {$p$}-subgroups of a group}, Adv. in Math., 28 (1978), pp.~101--128.

\bibitem{ar:S80}
{\sc J.-P. Serre}, {\em Trees}, Springer-Verlag, Berlin, 1980.
\newblock Translated from the French by John Stillwell.

\bibitem{ar:Spp}
{\sc R.~K. Skora}, {\em Deformation of length functions in groups}.
\newblock preprint.

\bibitem{ar:T70}
{\sc J.~Tits}, {\em Sur le groupe des automorphismes d'un arbre}, in Essays on
  topology and related topics (M\'emoires d\'edi\'es \`a Georges de Rham),
  Springer, New York, 1970, pp.~188--211.

\end{thebibliography}

\end{document}